\RequirePackage[leqno]{amsmath}
\documentclass[leqno]{amsart}
\usepackage{amsmath,mathtools}
\makeatletter
\usepackage[utf8]{inputenc}     
\usepackage{multicol}
\usepackage[inline]{enumitem} 
\usepackage{color, graphics}
\usepackage[demo]{graphicx}
\usepackage{efbox,graphicx}
\efboxsetup{linecolor=gray,linewidth=0.8pt}
\usepackage{floatrow}
\usepackage{float}
\usepackage{pgfplots}
\usepackage{mathrsfs}
\usepackage{latexsym, amsmath, amssymb}
\usepackage{caption}
\usepackage{subcaption}
\usepackage{pgfplots}
\usepackage{amsfonts}
\usepackage{pifont}
\usepackage{enumitem}
\setlist[itemize]{topsep=0pt,after=\vspace{1.5\baselineskip}}
\usepackage[colorlinks=true]{hyperref}
\usepackage{empheq}
\usepackage[T1]{fontenc}
\usepackage{tabularx}
\usepackage[leftcaption]{sidecap}
\sidecaptionvpos{figure}{m}
\usepackage{tikz}
\usepackage{caption}
\usepackage{tikz-cd}
\usetikzlibrary{decorations.markings}
\usetikzlibrary{decorations.pathmorphing}
\usepackage{tabularx,multirow,array} 
\usepackage[margin=0.5in]{geometry}
\usetikzlibrary{backgrounds,automata}

\pgfplotsset{compat=1.18}

\usepackage{xparse}
\ExplSyntaxOn
\int_new:N \g_tohi_const_int
\int_new:N \g_tohi_const_sub_int
\tl_new:N  \g_tohi_const_char_tl

\cs_new_protected:Nn \tohi_print_constant:nn 
{
  #1 \textsubscript {#2} 
} 

\NewDocumentCommand\resetconstants{m}
{
 \int_gincr:N \g_tohi_const_int
 \int_gzero:N \g_tohi_const_sub_int
 \tl_gset:Nn  \g_tohi_const_char_tl {#1}
}

\NewDocumentCommand\const{m}
{
  \tl_if_exist:cTF
   {
    c_tohi_const_\int_use:N\g_tohi_const_int _#1_tl
   }
   {
    \tl_use:c {c_tohi_const_\int_use:N\g_tohi_const_int _#1_tl }
   }
   {
    \int_gincr:N \g_tohi_const_sub_int
    \tl_const:cx {c_tohi_const_\int_use:N\g_tohi_const_int _#1_tl }
     { \exp_not:N\tohi_print_constant:nn {\g_tohi_const_char_tl }{\int_use:N \g_tohi_const_sub_int}}
    \tl_use:c {c_tohi_const_\int_use:N\g_tohi_const_int _#1_tl }
   }
}

\ExplSyntaxOff
\newcommand{\inlineitem}[1][]{%
\ifnum\enit@type=\tw@
    {\descriptionlabel{#1}}
  \hspace{\labelsep}%
\else
  \ifnum\enit@type=\z@
       \refstepcounter{\@listctr}\fi
    \quad\@itemlabel\hspace{\labelsep}%
\fi}
\DeclarePairedDelimiter\abs{\lvert}{\rvert}
\DeclarePairedDelimiter\norm{\lVert}{\rVert}
\DeclarePairedDelimiter\tonda{(}{)}
\DeclarePairedDelimiter\quadra{[}{]}
\DeclarePairedDelimiter\graffa{\{}{\}}
\newcommand{\into}{\int_\Omega}

\setlist[itemize]{noitemsep, topsep=0pt}

\def\R{\mathbb R} \def\N{\mathbb N}

\def\R{\mathbb R} \def\N{\mathbb N} 
\def\TM{T_{\textup{max}}} 
\def
\@cite
#1#2{[#1\if@tempswa, #2\fi]}
\makeatother
\newtheorem{theorem}{Theorem}[section]
\newtheorem{corollary}[theorem]{Corollary}

\newtheorem{lemma}[theorem]{Lemma}

\newtheorem{remark}{Remark}
\newtheorem{criterion}{Criterion}

\title[Boundedness in a class of Keller--Segel models with dissipative gradient terms] 
{
Dissipative gradient nonlinearities prevent $\delta$-formations in local and nonlocal attraction-repulsion chemotaxis models 
}
\author[Tongxing Li, Daniel Acosta-Soba, Alessandro Columbu, Giuseppe Viglialoro]{$^{\sharp}$Tongxing Li, $^{\flat}$Daniel Acosta-Soba, $^{\natural}$Alessandro Columbu and $^{\natural}$Giuseppe Viglialoro$^{\star}$}
\keywords{Chemotaxis, Attraction-repulsion, Nonlinear production, Gradient nonlinearities, Boundedness, Simulations. \\
\textit{$^\star$Corresponding author}:giuseppe.viglialoro@unica.it}
\begin{document}
\maketitle
{
\medskip
\centerline{$^{\sharp}$School of Control Science and Engineering}
\centerline{Shandong University}
\centerline{Jinan, Shandong, 250061 (P. R. China)}
\medskip
\centerline{$^{\flat}$Departamento de Matemáticas}
\centerline{Universidad de Cádiz, Puerto Real} \centerline{Campus Universitario Río San Pedro s/n, 11510. Cádiz (Spain)}
\medskip
\centerline{$^{\flat}$Department of Mathematics}
\centerline{University of Tennessee at Chattanooga, Chattanooga}
\centerline{615 McCallie Ave, 37403. TN (USA)}
\medskip
\centerline{$^{\natural}$Dipartimento di Matematica e Informatica}
\centerline{Universit\`{a} di Cagliari}
\centerline{Via Ospedale 72, 09124. Cagliari (Italy)}
\medskip
}
\bigskip
\resetconstants{c}
\begin{abstract}
We study a class of zero-flux attraction-repulsion chemotaxis models, characterized by nonlinearities laws  for the diffusion of the cell density $u$, the chemosensitivities and the production rates of the chemoattractant $v$ and the chemorepellent $w$. Additionally, a source involving also the gradient of $u$ is incorporated. The related mathematical formulation reads as:
\begin{equation}\label{problem_abstract}
\tag{$\Diamond$}
\begin{cases}
u_t= \nabla \cdot \tonda*{(u+1)^{m_1-1}\nabla u -\chi u(u+1)^{m_2-1}\nabla v + \xi u(u+1)^{m_3-1}\nabla w} +\lambda u^\rho -\mu u^k - c\abs{\nabla u}^\gamma  & \textrm{ in } \Omega \times (0,\TM),\\
\mathscr{P}(v,u)=\mathscr{Q}(w,u)=0 & \textrm{ in } \Omega \times (0,\TM).
\end{cases}
\end{equation}
Herein $\Omega$ is a bounded and smooth domain of $\R^n$, for $n\in \N$, $\TM\in (0,\infty], \chi,\xi,\rho,k,\gamma,\lambda,\mu,c$ are proper positive numbers and $m_1,m_2,m_3\in \R$. 
Depending on the expressions of $\mathscr{P}(v,u)$ and $\mathscr{Q}(w,u)$, we will deal with issues tied to \textit{local} or \textit{nonlocal} models. 

Our overall study  touches different aspects: we address questions connected to local well-posedness, we derive sufficient conditions so to ensure boundedness of solutions and finally we develop numerical simulations giving insights on the evolution of the system.

In the specific, for $\tau\in\{0,1\}$, and $f_1(u)$ and $f_2(u)$ reasonably regular functions generalizing, respectively, the prototypes $f_1(u)= u^\alpha$ and $f_2(u)= u^\beta$, for some $\alpha,\beta>0$ and all $u\geq 0$, by defining $\Theta_\tau=\max\left\{1,\frac{n}{n+1}(m_2+\alpha),\tau\frac{n}{n+1}(m_3+\beta)\right\}$ this will be established for model \eqref{problem_abstract}:
\begin{enumerate}[label=(\Roman*),ref=\Roman*]
\item \label{LocalWellItem} For $\mathscr{P}(v,u)=P^\tau_\alpha(v,u):=\tau v_t-\Delta v+v-f_1(u)$ and $\mathscr{Q}(w,u)=Q^\tau_\beta(w,u):=\tau w_t-\Delta w+w-f_2(u)$ (i.e., the local case), and for $\mathscr{P}(v,u)=P_\alpha(v,u):=\Delta v-\frac{1}{|\Omega|}\int_\Omega f_1(u)+f_1(u)$ and $\mathscr{Q}(w,u)=Q_\beta(v,u):=\Delta w-\frac{1}{|\Omega|}\int_\Omega f_2(u)+f_2(u)$  (i.e., the nonlocal case), $\rho,k,\gamma\geq 1$,  whenever $u_0,\tau v_0, \tau w_0$ belong to specific Hölder spaces, the system \eqref{problem_abstract} admits for some $0<T<\TM$ a unique classical solution in $\bar\Omega \times [0,T]$, with some Hölder regularity in space and time.
\item \label{BoundedItem} For  $\mathscr{P}(v,u)=P^\tau_\alpha(v,u)$, $\mathscr{Q}(w,u)=Q^\tau_\beta(w,u)$ and $u_0,\tau v_0, \tau w_0$ as in \eqref{LocalWellItem}, for $1\leq \rho<k$ and  $\Theta_\tau<\gamma\leq 2$, the system  is globally solvable (i.e., $\TM=\infty$) and the solution is bounded.
\item For $\mathscr{P}(v,u)=P_\alpha(v,u)$ and $\mathscr{Q}(w,u)=Q_\beta(w,u)$, and $u_0$ as in \eqref{LocalWellItem}, for $\tau=0$, $1\leq \rho<k$ and  $\Theta_\tau<\gamma\leq 2$ the same conclusion as in \eqref{BoundedItem} applies.
\item Numerical simulations show that the suppression of the condition $\Theta_\tau<\gamma\leq 2$ may lead to only local solutions (i.e., $\TM$ finite) which precisely blow-up at $\TM$.
\end{enumerate}
This paper, inter alia, generalizes the result in  \cite{IshidaLankeitVigliloro-Gradient}, where the linear and only attraction version  of model \eqref{problem_abstract} is addressed.
\end{abstract}
\section{Introduction, motivations, state of the art and aim of the research}\label{Intro}
\subsection{The continuity equation with dissipative gradient terms \texorpdfstring{$\boldsymbol{h=h(u,\nabla u)}$}{h(u)} in biological models}
The \textit{continuity equation} 
\begin{equation}\label{ContinuityEq}
    u_t+\nabla \cdot F=h
\end{equation}
nonnegative physical quantity, denoted as $u=u(x,t)$, at position $x$ and time $t>0$. Within this equation, the flux $F$ depicts the flow of this quantity, while $h$ denotes an external influence making that $u$ may increase and/or decrease over time.

The main emphasis of this paper lies in considering chemotaxis models tied to \eqref{ContinuityEq} in which the source $h$ presents effects of dissipative type involving gradient terms. Indeed, gradient-dependent sources may influence the dynamics of the  quantity obeying the continuity equation. Confining our attention to biological applications, according to \cite{Souplet_Gradient}, a biological species with density $u$ and occupying a certain habitat evolves in time by displacement, birth/reproduction and death. In particular, the births are described by a superlinear power of such a distribution, the natural deaths by a linear one and the accidental deaths by a function of its gradient; this leads (through the choice $F=F(u)=-\nabla u$ and $h=h(u,\nabla u)=u^\rho-u- |\nabla u|^\gamma$ in \eqref{ContinuityEq}) to $u_t=\Delta u +u^\rho-u-|\nabla u|^\gamma$, with $\rho,\gamma>1$.  

The study of gradient terms in parabolic equations, and their impact on the possibility of blow-up, in the sense that $u$ tends to a coalescence and ceases existing at some finite time (we will be more clear below),  is a classical topic;
we mention \cite{chipot_weissler,kawohl_peletier,fila} and \cite[Sect. IV]{QS}, where the relations between exponents $\rho,\gamma$  in the Dirichlet problem for a semilinear heat equation of the form $u_t=\Delta u +u^\rho-|\nabla u|^{\gamma}$  preventing or implying  blow-up are investigated.  
(We refer also to \cite{SoupletRecent2011,Souplet2005Book} for surveys and results on the topic.)



\subsection{A review on chemotaxis mechanisms with \texorpdfstring{$\boldsymbol{h\equiv 0}$}{h0},  \texorpdfstring{$\boldsymbol{h=h(u)}$}{h(u)} and \texorpdfstring{$\boldsymbol{h=h(u,\nabla u)}$}{h(u)}}
This study focuses on the analysis of  equation \eqref{ContinuityEq} within the framework of chemotaxis phenomena observed in biological populations. These mechanisms involve organisms or entities altering their trajectory in response to one or more chemical signals. Specifically, our interest lies in the dynamics of a certain cell density $u=u(x,t)$, for which the flux $F$ comprises a smooth diffusive component (connected to $-\nabla u$) and a double-action field (related to $u\nabla v-u\nabla w$), arising from two chemical signals: a chemoattractant $v=v(x,t)$, which tends to aggregate cells, and a chemorepellent $w=w(x,t)$, which disperses them. As to the source $h$, oppositely to the case of a single species, as far as we know, the actual literature on mathematical analysis concerning chemotaxis models with gradient depending sources is far to be totally discussed and understood (we will specify it below). For this reason, herein we incorporate an external \textit{gradient-dependent source} that influences the cell density and may exert on it both increasing and decreasing effects.  
\subsubsection{The attraction-repulsion mechanism}
Oppositely to what now said, the knowledge concerning chemotaxis mechanisms without sources $h$ or with $h=h(u)$, only, is rather deep. In the specific, wanting to give a presentation in terms of the continuity equation \eqref{ContinuityEq}, it is meaningful (and convenient) to define, for positive parameters $\chi$, $\xi$, $\lambda$, $\mu$, and general real $m_1$, $m_2$, $m_3$, along with $k>1$, the fluxes $F=F_{m_1,m_2,m_3}=F_{m_1,m_2,m_3}(u,v,w)$, $G_{m_1}=G_{m_1}(u)$, $H_{m_2}=H_{m_2}(u,v)$, $I_{m_3}=I_{m_3}(u,w)$, and the source $h=h_k=h_k(u)$:
\begin{equation}\label{FluxAndSource}
\begin{cases}
F=F_{m_1,m_2,m_3}=-(u+1)^{m_1-1}\nabla u +\chi u(u+1)^{m_2-1}\nabla v-\xi u(u + 1)^{m_3-1}\nabla w=:G_{m_1} +H_{m_2}+I_{m_3},\\
h= h_k=\lambda u -\mu u^k.
\end{cases}
\end{equation}
Here, the diffusion $-(u+1)^{m_1-1}\nabla u$ becomes more significant for higher values of $m_1$, while the aggregation/repulsion effects, $\chi u(u+1)^{m_2-1}\nabla v$ and $-\xi u(u + 1)^{m_3-1}\nabla w$, increase with larger sizes of $\chi$, $\xi$, $m_2$, and $m_3$. Additionally, the cell density may increase at a rate of $\lambda u$ and decrease at a rate of $-\mu u^k$. The flux is influenced by the signals $v$ and $w$,
which evolve according to two differential equations that
we currently denote as $\mathscr{P}(v)=0$ and $\mathscr{Q}(w)=0$ (to be specified later).
One of these equations describes the behavior of the chemoattractant $v$, the other for the chemorepellent $w$, and both have to be  coupled with the continuity equation. Furthermore, in order to make the analysis complete, a boundary condition and  the initial configurations for the cell and chemical densities have to be assigned. Essentially, for 
insulated domains, with homogeneous Neumann or zero flux boundary conditions,  we are concerned (taking equations \eqref{ContinuityEq} and \eqref{FluxAndSource} in mind) with the following initial boundary value problem:
\begin{equation}\label{problemAstratto}
\begin{cases}
u_t+\nabla \cdot F_{m_1,m_2,m_3}=h_k & \text{ in } \Omega \times (0,T_{\text{max}}), \\
\mathscr{P}(v,u)=\mathscr{Q}(w,u)=0 & \text{ in } \Omega \times (0,T_{\text{max}}),\\
u_0(x)=u(x,0)\geq 0; v_0(x)=v(x,0)\geq 0; w_0(x)=w(x,0) \geq0 & x \in \bar\Omega,\\
u_{\nu}=v_{\nu}=w_{\nu}=0 & \text{ on } \partial \Omega \times (0,T_{\text{max}}).
\end{cases}
\end{equation}
This problem is formulated in a bounded and smooth domain $\Omega$ of $\mathbb{R}^n$, where $n\geq 1$, and $u_{\nu}$ (and similarly for $v_{\nu}$ and $w_{\nu}$) indicates the outward normal derivative of $u$ on $\partial \Omega$, the boundary of $\Omega$. Moreover, $T_{\textrm{max}}\in (0,\infty]$ identifies the maximum time up to which solutions to the system can be extended.

The discussion above finds its roots in the well-established Keller--Segel models, which idealize chemotaxis phenomena and have been of interest to the mathematical community for the past 50 years, as evidenced by the seminal papers \cite{K-S-1970,Keller-1971-MC,Keller-1971-TBC}.

On the other hand, since we will deal with \textit{local} and \textit{nonlocal} biological models,  we summarize some known results connected to such situations.
\subsubsection{Boundedness vs. blow-up: the local case}
\textit{Local} dynamics relate to nearby activities, and quantities involved in the related models are influenced directly only by factors in their immediate surroundings. 

Specifically, concerning chemotaxis models with a single proliferation signal, considering \eqref{FluxAndSource}, problem \eqref{problemAstratto} represents a more general combination of an aggregative signal-production mechanism
\begin{equation}\label{problemOriginalKS}
u_t= -\nabla \cdot \left(G_1+H_1\right)= \Delta u- \chi \nabla \cdot (u \nabla v)\quad \textrm{and} \quad
\mathscr{P}(v,u)=P^\tau_1(v,u)=\tau v_t-\Delta v+v-u=0, \quad \textrm{ in } \Omega \times (0,T_{\textrm{max}}),
\end{equation}
and a repulsive signal-production counterpart 
\begin{equation}\label{problemOriginalKSCosnumption}
u_t= -\nabla \cdot \left(G_1+I_1\right)= \Delta u+ \xi \nabla \cdot (u \nabla w) \quad \text{and} \quad
\mathscr{Q}(w,u)=Q^\tau_1(w,u)=\tau w_t-\Delta w+w -u=0, \quad \text{ in } \Omega \times (0,T_{\text{max}}).
\end{equation}
Here, $\tau\in\{0,1\}$ distinguishes between a stationary and an evolutive equation for the chemical. These models exhibit linear diffusion (i.e., $\nabla \cdot \nabla u=\Delta u$) and linear production rates: specifically, $v$ and $w$ are linearly produced by the cells themselves. The opposite phenomenon concerns the situation where the particles' density consumes the chemical (models with absorption will only be briefly mentioned); formally, in  $P^1_1(v,u)$ the term $v-u$ is replaced by $uv$ (or in $Q^1_1(w,u)$ by $uw$). 

As for problem \eqref{problemOriginalKS}, since the attractive signal $v$ increases with $u$, the natural spreading process of the cells' density could be interrupted, leading to the formation of highly concentrated spikes (\textit{chemotactic collapse} or \textit{blow-up at finite time}). This instability generally is connected  to the size of the chemosensitivity $\chi$, the initial mass of the cell distribution, denoted as $m=\int_\Omega u_0(x)dx$, and the space dimension $n$. In this regard, for details we refer to \cite{HerreroVelazquez,JaLu,Nagai,WinklAggre}, where analyses concerning the existence and properties of global, uniformly bounded, or local solutions to similar models as \eqref{problemOriginalKS} can be found.

On the contrary, for chemotaxis models with nonlinear segregation, such as the ones under our focus, uniform-in-time boundedness of all solutions to problem \eqref{problemOriginalKS} with the choice $\mathscr{P}(v,u)=P^1_\alpha(v,u)=v_t-\Delta v+v-u^\alpha$, with $0<\alpha<\frac{2}{n}$ ($n\geq1$), is shown in  \cite{LiuTaoFullyParNonlinearProd}.

Regarding the literature concerning problem \eqref{problemOriginalKSCosnumption}, it appears to be poor (see, for instance, \cite{Mock74SIAM,Mock75JMAA} for analyses on similar contexts). Particularly, no results on  blow-up scenarios are available, and this is  intuitive due to the repulsive nature of the phenomenon.

In contrast, the understanding level for attraction-repulsion chemotaxis problems involving both \eqref{problemOriginalKS} and \eqref{problemOriginalKSCosnumption} is considerably rich. Specifically, in the framework of linear diffusion and sensitivities version of model \eqref{problemAstratto}, where $F=F_{1,1,1}$ and $\mathscr{P}(v,u)=P^\tau_\alpha(v,u)=\tau v_t-\Delta v+b v -a u^\alpha$ and $\mathscr{Q}(w,u)=Q^\tau_\beta(w,u)=\tau w_t- \Delta w + d w - c u^\beta$, with $a,b,c,d,\alpha,\beta>0$, and equipped with regular initial data $u_0(x),\tau v_0(x),\tau w_0(x)\geq 0$, several outcomes can be summarized. For stationary contexts, and in the absence of logistics ($h_k\equiv 0$), when linear growths of the chemoattractant and the chemorepellent are considered, and with elliptic equations for the chemicals (i.e., when $P^0_1(v,u)=Q^0_1(w,u)=0$), the value $\Theta:=\chi a-\xi c$ quantifies the difference between the attraction and repulsion impacts. It is observed that whenever $\Theta<0$ (repulsion-dominated regime), all solutions to the model are globally bounded in any dimension. Conversely, for $\Theta>0$ (attraction-dominated regime) and $n=2$, unbounded solutions can be detected. For more general expressions of the proliferation laws, described by the equations $P^0_\alpha(v,u)=Q^0_\beta(w,u)=0$,  in \cite{ViglialoroMatNacAttr-Repul} some interplay between $\alpha$ and $\beta$ (along with certain technical conditions on $\xi$ and $u_0$) are established to ensure the globality and boundedness of classical solutions. Additionally, \cite{ChiyoYokotaBlow-UpAttRe} provides blow-up results within the framework of nonlinear attraction-repulsion models with logistics, akin to those formulated in \eqref{problemAstratto} with $F=F_{m_1,m_2,m_3}$ and $h=h_k$, and with linear segregation for the stimuli, i.e., with equations for $v$ and $w$ reading as $P_1^0(v)=Q_1^0(w)=0$.

Focusing on evolutive equations for the chemoattractant and chemorepellent, specifically $P^1_1(v,u)=Q^1_1(w,u)=0$, \cite{TaoWanM3ASAttrRep} establishes that in two-dimensional domains with sufficiently smooth initial data, globally bounded solutions persist over time whenever
$$\Theta<0 \quad \textrm{and} \quad b=d  \quad \textrm{or} \quad  \Theta<0 \quad \textrm{and} \quad  -\frac{\chi^2\alpha^2(b-d)^2}{2\Theta {b}^2 C} \int_\Omega u_0(x)dx\leq 1, \quad \textrm{for some }\; C>0.$$

Regarding blow-up results, we are only aware of \cite{LankeitJMAA-BlowupParabolicoAttr-Rep}, where unbounded solutions in three-dimensional domains are constructed.

When $h=h_k\not\equiv 0$, for both linear and nonlinear production scenarios, and stationary or evolutive equations (formally, $P^{\tau}_\alpha(v,u)=Q^{\tau}_\beta(w,u)=0$), criteria concerning boundedness, long-term behavior, and blow-up issues for related solutions are investigated in \cite{LiangEtAlAtt-RepNonLinProdLogist-2020,XinluEtAl2022-Asymp-AttRepNonlinProd,ChiyoMarrasTanakaYokota2021,GuoqiangBin-2022-3DAttRep}.
\subsubsection{Boundedness vs. blow-up: the nonlocal case}
In contrast to the local case, \textit{nonlocal} dynamics involve quantities/events that affect distant areas, beyond their immediate vicinity. 
As to such models associated with \eqref{problemAstratto}, the nonlocality action is generally connected to the production laws;  specifically $\mathscr{P}(v,u)=P_\alpha(v,u):=\Delta v-\frac{1}{|\Omega|}\int_\Omega u^\alpha +u^\alpha=0$, and $\mathscr{Q}(w,u)=Q_\beta(w,u):=\Delta w-\frac{1}{|\Omega|}\int_\Omega u^\beta +u^\beta$.
In this scenario, the unknown $v=v(x,t)$  represents the chemoattractant \textit{deviation}, and not the distribution itself as in the local models. In the specific we have that $v$ changes sign, and its  mean value is zero, exactly as it will be imposed below.  (Naturally, the chemical repellent $w$ behaves similarly; we do not make use of a different symbolism for the deviations, since from the context it is clear to which quantity we are referring to.)  Particularly, following the convention in \eqref{FluxAndSource}, we now highlight some results on nonlocal phenomena  inspiring our investigation. To this end, we discuss two specific cases: the sole attraction version
\begin{equation}\label{OnlyATTractionproblem}
u_t= -\nabla \cdot \left(G_{m_1}+H_{m_2}\right) +h_k \quad \textrm{and} \quad
\mathscr{P}(v,u)=P_\alpha(v,u)=0 \quad \textrm{ in } \Omega \times (0,T_{\textrm{max}}),
\end{equation}
and the attraction-repulsion model
\begin{equation}\label{ATTract-Repulsionproblem}
u_t= -\nabla \cdot F_{m_1,m_2,m_3} +h_k \quad \textrm{and} \quad
\mathscr{P}(v,u)=P_\alpha(v,u)=0 \textrm{ and } \mathscr{Q}(w,u)=Q_\beta(w,u)=0 \quad \textrm{ in } \Omega \times (0,T_{\textrm{max}}).
\end{equation}
Regarding the linear diffusion and sensitivity version of model \eqref{OnlyATTractionproblem}, i.e., with flux $G_{1}+H_{1}$, it has been established that if $h_k\equiv 0$, solutions remain bounded for any $n\geq 1$ and $0<\alpha<\frac{2}{n}$. Conversely, for $\alpha>\frac{2}{n}$, blow-up phenomena may occur (refer to \cite{WinklerNoNLinearanalysisSublinearProduction} for further details). Introducing damping logistic terms into the mechanism, $\delta$-formations at finite time have as well been observed for the linear flux $G_1+H_1$ under certain sub-quadratic growth conditions of $h_k$, specifically for $h_k$ with $1<k<\frac{n(\alpha+1)}{n+2}<2$ (see \cite{yi2021blow}). Additionally, for the linear production situation, corresponding to $P_1(v)=0$, unbounded solutions have been detected in \cite{FuestCriticalNoDEA} even for quadratic sources $h=h_2$, provided $n\geq 5$ and $\mu\in \left(0,\frac{n-4}{n}\right)$.

In nonlinear models without damping logistic effects (i.e., with a general flux $G_{m_1}+H_{m_2}$ and $h_k\equiv 0$) and linear production (i.e., $P_1(v)=0$), among other findings established in \cite{winkler2010boundedness} we indicate that  for $m_1\leq 1$, $m_2>0$, and $m_2>m_1+\frac{2}{n}-1$, scenarios with unbounded solutions at certain finite times $T_{\textrm{max}}$ can be encountered. For related questions concerning estimates of $\TM$, see also \cite{MarrasNishinoViglialoro}. Additionally, some results have been further explored in \cite{TANAKA2021} for  $h_k\not\equiv 0$ and within nonlinear segregation contexts ($P_\alpha(v,u)=0$). Particularly, for $m_1\in \R$, $m_2>0$, blow-up phenomena are observed to occur if $ m_2+\alpha > \max\left\{m_1+\frac{2}{n}k, k\right\}$, whenever $m_1 \geq 0$, or $m_2+\alpha > \max\left\{\frac{2}{n}k, k\right\}$, provided $m_1 <0$.

Conversely, moving our attention to  attraction-repulsion models, in \cite{LiuLiBlowUpAttr-Rep} it is proved that if $P_\alpha(v,u)=Q_\beta(w,u)=0$, the conditions $\alpha>\frac{2}{n}$ and $\alpha>\beta$ guarantee the existence of unbounded solutions to \eqref{ATTract-Repulsionproblem} for the linear flux $F=F_{1,1,1}$, without logistic effects ($h_k\equiv 0$). However, identifying gathering mechanisms for the nonlinear scenario  and in the presence of damping logistics is more intricate and, to our knowledge, there are only two references addressing this issue. On the one hand, \cite{wang2023blow} focuses on $F=F_{m_1,1,1}$ with $m_1\in \R$. On the other hand, more recently, the existence of blow-up solutions for a broader range of fluxes, specifically for $F=F_{m_1,m_2,m_2}$ with $m_1\in \R$ and any $m_2=m_3>0$, have been detected in  \cite{ColFraVig-StudApplMat}.

\subsubsection{Gradient-dependent sources \texorpdfstring{$h=h(u,\nabla u)$}{h(u)}}
As above anticipated, in chemotaxis models gradient nonlinearities involved in external sources seem to have been poorly studied so far; we can only mention  \cite{ViglialoroDifferentialIntegralEquations} and  \cite{IshidaLankeitVigliloro-Gradient} (both framed in local scenarios), where  real biological interpretations also connected to  ecological  mechanisms are given. In particular, by introducing the symbol $h_{\rho,k,\gamma}=\lambda u^\rho-\mu u^k-c|\nabla u|^\gamma$ to describe a source with dissipative gradient terms for the model
\begin{equation}\label{OnlyATTractionproblemBis}
u_t= -\nabla \cdot \left(G_{m_1}+H_{m_2}\right) +h_{\rho,k,\gamma} \quad \textrm{and} \quad
P_\alpha^\tau(v,u)=0 \quad \textrm{ in } \Omega \times (0,T_{\textrm{max}}),
\end{equation}
 lower bounds of the maximal existence time of given solutions are derived in two- and three-dimensional settings, for some $\vartheta>2$, $p>q>2$, whenever $m_1=m_2=1$, $P_\alpha^1(v,u)$, for $\alpha=\vartheta \frac{p-1}{2}$ and $h_{\rho,k,\gamma}$, with $k=\rho=p, \gamma=q$ (\cite{ViglialoroDifferentialIntegralEquations}). Nevertheless, \cite{IshidaLankeitVigliloro-Gradient} deals with the existence and boundedness of classical solutions to problem \eqref{OnlyATTractionproblemBis} with $m_1=m_2=1$, and $P^0_\alpha (v)$ and  $P^\tau_1(v)$; in this sense, such boundedness of classical solutions is achieved in any dimension $n$, for any initial data, whenever $k>\rho$ and $\frac{2n}{n+1}<\gamma\leq 2.$ 
\subsection{Aim of the investigation}
Motivated by the discussed considerations indicating that generalized logistic sources may still permit blow-up even in attraction-repulsion models reading as
$$
u_t+ \nabla \cdot F_{m_1,m_2,m_3} =h_k,
$$
in this investigation
\textit{we aim at analyzing exactly these same unstable systems but perturbed by gradient-dependent sources, i.e. models mechanisms of the type 
$$
u_t+\nabla \cdot F_{m_1,m_2,m_3} =h_{\rho,k,\gamma},
$$
and at deriving sufficient conditions involving
the accidental deaths term (i.e., $\gamma$) such that blow-up is prevented.
}

In the specific 
\begin{itemize}
    \item [$\diamond$] we will focus on the theoretical analysis dealing with existence and boundedness of classical solutions to fully nonlinear attraction-repulsion models of local- and nonlocal-type, so extending that already developed in \cite{IshidaLankeitVigliloro-Gradient} for the only attraction, linear and local scenario;
    \item [$\diamond$] we will present some numerical simulation in order to have some hints on whether even weak damping gradient terms may lead to $\delta$-formations. 
\end{itemize}
\section{Presentation of the models and of the main results}
Accordingly to the above preparations, we herein are interested to address questions connected to local-well posedness and global solvability of the following \emph{local} and \emph{nonlocal} attraction-repulsion chemotaxis models with gradient-dependent sources: 
\begin{equation}\label{localproblem}
\tag{$\mathcal{L}_\tau$}
\begin{dcases}
u_t= \nabla \cdot \tonda*{(u+1)^{m_1-1}\nabla u -\chi u(u+1)^{m_2-1}\nabla v + \xi u(u+1)^{m_3-1}\nabla w} +\lambda u^\rho -\mu u^k - c\abs{\nabla u}^\gamma & \text{in } \Omega \times (0,\TM),\\
\tau v_t= \Delta v - v  +f_1(u) & \text{in } \Omega \times (0,\TM),\\
\tau w_t= \Delta w -  w + f_2(u) & \text{in } \Omega \times (0,\TM),\\
u_{\nu}=v_{\nu}=w_{\nu}=0 & \text{on } \partial \Omega \times (0,\TM),\\
u(x,0)=u_0(x), \tau v(x,0)=\tau v_0(x), \tau w(x,0)=\tau w_0(x)  & x \in \bar\Omega,
\end{dcases}
\end{equation}
and
\begin{equation}\label{nonlocalproblem}
\tag{$\mathcal{NL}$}
\begin{dcases}
u_t= \nabla \cdot \tonda*{(u+1)^{m_1-1}\nabla u -\chi u(u+1)^{m_2-1}\nabla v + \xi u(u+1)^{m_3-1}\nabla w}+\lambda u^\rho -\mu u^k - c\abs{\nabla u}^\gamma & \text{in } \Omega \times (0,\TM),\\
0= \Delta v - \frac{1}{\abs{\Omega}}\into f_1(u)  + f_1(u) & \text{in } \Omega \times (0,\TM),\\
0= \Delta w - \frac{1}{\abs{\Omega}}\into f_2(u)  + f_2(u) & \text{in } \Omega \times (0,\TM),\\
u_{\nu}=v_{\nu}=w_{\nu}=0 & \text{on } \partial \Omega \times (0,\TM),\\
u(x,0)=u_0(x) & x \in \bar\Omega,\\
\int_\Omega v(x,t)\,dx=\int_\Omega w(x,t)\,dx=0 & \text{for all } t\in  (0,\TM).
\end{dcases}
\end{equation}
Here, we assume  
\begin{equation}\label{reglocal}
\begin{cases}
   f_1,f_2: [0,\infty) \rightarrow \R^+, \textrm{ with } f_1,f_2\in C^1([0,\infty)), 
   \\ 
    \textrm{ for some } \delta\in(0,1) \textrm{ and } n\in \N, \Omega \subset \R^n \textrm{ is a bounded domain of class } C^{2+\delta}, \textrm{ with boundary } \partial \Omega,\\
    u_0, v_0, w_0: \bar{\Omega}  \rightarrow \R^+ , \textrm{ with } u_0,  v_0,   w_0 \in C_\nu^{2+\delta}(\bar\Omega)=\{\phi \in C^{2+\delta}(\bar{\Omega}): \phi_\nu=0 \textrm{ on }\partial \Omega\}. 
   \end{cases}
\end{equation}
Moreover, we might require that for proper $\alpha,\beta,k_1,k_3 \geq k_2>0,$
\begin{equation}\label{disf}
    f_1(s)\leq k_1 (s+1)^\alpha \quad \text{and} \quad  k_2(s+1)^\beta \leq f_2(s)\leq k_3 (s+1)^\beta.
\end{equation}
Finally, for $\tau\in\graffa{0,1}$, $m_1,m_2,m_3\in\R$ and $\alpha,\beta,\chi,\xi,\lambda,\mu,c>0$, $1\leq\rho< k$ let 
\begin{equation}\label{condgamma}  
\max{\graffa*{1,\frac{n}{n+1}(m_2+\alpha),\tau\frac{n}{n+1}(m_3+\beta)}}<\gamma\leq2. \tag{$\mathcal{G}_\tau$}
\end{equation} 
We will prove the following results.
\begin{theorem}\label{theoremExistence}
For $\tau\in \{0,1\}$, let $\Omega$, $f_1,f_2, u_0, \tau v_0, \tau w_0$ comply with hypotheses in \eqref{reglocal}.  Additionally, let $\chi,\xi,\lambda,\mu,c>0$, and $m_1,m_2,m_3\in \R$, $k,\rho,\gamma \geq 1$.  Then there exist $T>0$ and a unique triple of functions $(u,v,w)$, with 
\begin{equation*}
(u,v,w)\in C^{2+\delta,1+\frac{\delta}{2}}( \Bar{\Omega} \times [0, T])\times C^{2+\delta,\tau+\frac{\delta}{2}}( \Bar{\Omega} \times [0, T])\times  C^{2+\delta,\tau+\frac{\delta}{2}}( \Bar{\Omega} \times [0, T]),
\end{equation*}
solving problems  \eqref{localproblem} and \eqref{nonlocalproblem}.
\end{theorem}
\begin{theorem}\label{theoremlocal}
For $\tau\in \{0,1\}$, let $\Omega$, $f_1,f_2, u_0,\tau v_0,\tau w_0$ comply with hypotheses in \eqref{reglocal} and \eqref{disf}.  
Additionally, let $\chi,\xi,\lambda,\mu,c>0$, and $m_1,m_2,m_3\in \R$, $k,\rho,\gamma \geq 1$
such that  condition \eqref{condgamma} hold true. Then problems \eqref{localproblem} and \eqref{nonlocalproblem} admit a unique solution
\begin{equation*}
(u,v,w)\in C^{2+\delta,1+\frac{\delta}{2}}( \Bar{\Omega} \times [0, \infty))\times C^{2+\delta,\tau+\frac{\delta}{2}}( \Bar{\Omega} \times [0, \infty))\times  C^{2+\delta,\tau+\frac{\delta}{2}}( \Bar{\Omega} \times [0, \infty)) 
\end{equation*}
such that $u,v,w \in L^\infty(\Omega \times (0,\infty)).$ Additionally $u,v,w\geq 0$ on $\bar\Omega \times [0,\infty)$, and $u\geq 0$ on $\bar\Omega \times [0,\infty)$,  for  \eqref{localproblem}, and \eqref{nonlocalproblem}, respectively. 
\end{theorem}
Let us give some comments.
\begin{remark}[Comparison with \cite{IshidaLankeitVigliloro-Gradient}]\label{RemarkDeviation}
Let us analyze our research in the spirit of \cite{IshidaLankeitVigliloro-Gradient}.
\begin{itemize}
\item [$\circ$] For $w\equiv 0$, models \eqref{localproblem} generalize those studied in \cite{IshidaLankeitVigliloro-Gradient}. In particular, for $m_1=m_2=\alpha=1$, problem \eqref{localproblem} is reduced to \cite[(2)]{IshidaLankeitVigliloro-Gradient} and condition \eqref{condgamma} reads $\frac{2n}{n+1}<\gamma \leq 2$; henceforth, \cite[Theorem 1.2]{IshidaLankeitVigliloro-Gradient} is a particular case of Theorem \ref{theoremlocal}.
\item [$\circ$] Even though the local-well posedness is established in \cite[Theorem 1.1]{IshidaLankeitVigliloro-Gradient}, for the fully nonlinear case some adjustments are required; additionally, conversely to \cite[Theorem 1.1]{IshidaLankeitVigliloro-Gradient}, Theorem \ref{theoremExistence} gives all the details concerning the uniqueness of classical solutions.
\item [$\circ$] \cite[Remark 1.5]{IshidaLankeitVigliloro-Gradient} leaves open the question  whether for some small $\gamma\in[1,\frac{2n}{n+1}]$ there may exist
unbounded solutions to problem \eqref{localproblem}, naturally for $w\equiv 0, m_1=m_2=\alpha=1$. Through some numerical simulations, we herein give evidences supporting some possible scenarios in such a direction. 
\item [$\circ$] In \cite{IshidaLankeitVigliloro-Gradient} the model involves a larger class of  sources than those herein utilized, and whose general expression is $h(u,\nabla u)=f(u)-g(\nabla u)$, being $f$ and $g$ function with certain growth properties. 
\end{itemize}
\end{remark}
\section{Necessary tools: parameters and properties of solutions to differential equations}
We will make use of these results, which we present without any extra comment. We also underline that 
\begin{itemize} 
\item [$\circ$] Without any mention, all the constants $c_i$, $i=1,2,\ldots$ are assumed to be positive;
\item [$\circ$] With $\varepsilon$ we indicate an \textit{arbitrary} positive real number, and the multiplication by another constant and the sum with  another \textit{arbitrary} constant are not performed, and the final result is for commodity as well labeled with $\varepsilon$. 
\end{itemize}
\begin{lemma}\label{BoundsInequalityLemmaTecnicoYoung}
For all $A,B \geq 0$ and $s>0$ we have
\begin{equation}\label{disab}
(A+B)^s \leq \max\{1,2^{s-1}\}(A^s+B^s).
\end{equation}
Moreover, if $\const{123}+\const{321}<1$, then for all  
$\varepsilon>0$ there is $c_3=c_3(\varepsilon)>0$ such that
\begin{equation}\label{LemmaEsponenti}
A^{\const{123}}B^{\const{321}} \leq \varepsilon A+B+\const{d}.
\end{equation}
\begin{proof}
The proofs are available in \cite[Theorem 1]{Jameson_2014Inequality}  and \cite[Lemma 4.3]{frassuviglialoro1}.
\end{proof}
\end{lemma}
\begin{lemma}\label{regularity}
Let $\Omega\subset\R^n$ satisfy condition in \eqref{reglocal}. 
\begin{itemize}
    \item [$\diamond$] {\bf\emph{Elliptic regularity}:}  If $\psi\in C^{\delta}(\bar\Omega)$, then the solution $z\in C^{2+\delta}(\bar\Omega)$  of
\begin{equation*}
\begin{cases}
    -\Delta z + z=\psi  & \text{in } \Omega,\\
    z_{\nu}=0 & \text{on } \partial\Omega,
\end{cases}
\end{equation*}
has the property that for every $q>1$ 
and any $\varepsilon>0$, there is $C_\varepsilon>0$ such that
\begin{equation}\label{ellipticreg}
    \into z^q\leq \varepsilon \into \psi^q + C_\varepsilon\tonda*{\into \psi}^q.
\end{equation}
 \item [$\diamond$] {\bf\emph{Parabolic regularity}:} For every $q>n$, $T\in(0,\infty]$, $g\in L^q([0,T);L^q(\Omega))$ and $z_0\in W^{2,q}(\Omega)$, with $\partial_{\nu} z_0 = 0$ on $\partial\Omega$, every solution $z\in W^{1,q}_{loc}([0,T);L^q(\Omega))\cap L^q_{loc}([0,T);W^{2,q}(\Omega))$ of 
\begin{equation*}
\begin{cases}
    z_t= \Delta z - z  +g & \text{in } \Omega \times (0,T),\\
    \partial_{\nu} z = 0 & \text{on } \partial\Omega \times (0,T),\\
    z(\cdot,0)=z_0 & \text{on } \Omega,
\end{cases}
\end{equation*}
satisfies for some $C_P>0$
\begin{equation}\label{tau1}
\int_0^t e^s \into \abs*{\Delta z(\cdot,s)}^q\,ds \leq C_{P}\quadra*{1+\int_0^t e^s \into\abs{g(\cdot,s)}^q\,ds} \quad \text{for all } t\in(0,T).
\end{equation}
Additionally, if  $g\in L^\infty([0,T);L^q(\Omega))$ then
\begin{equation}\label{tau1extension}
    \nabla z \in L^\infty((0,T);L^\infty(\Omega)).
\end{equation}
\end{itemize}
\begin{proof} A detailed proof of \eqref{ellipticreg} is available in \cite[Lemma 3.1]{ViglialoroMatNacAttr-Repul}.
For \eqref{tau1} we refer to \cite{IshidaLankeitVigliloro-Gradient}, whereas for \eqref{tau1extension} we can invoke \cite[Lemma 4.1]{HorstWink} and the embedding $W^{1,\frac{n q}{(n-q)_+}}(\Omega)\subset L^\infty(\Omega)$, valid for $q>n$.
\end{proof}
\end{lemma}
\begin{lemma}\label{LemmaODI-Comparison}
Let\/ $T>0$ and $\phi:(0,T)\times \R^+_0\rightarrow \R$. If $0\leq y\in C^0([0,T))\cap  C^1((0,T))$ is such that 
\begin{equation*}
y'\leq \phi(t,y)\quad \textrm{for all } t \in (0,T), 
\end{equation*}
and there is $y_1>0$ with the property that whenever $y>y_1$ for some $t\in (0,T)$ one has that $\phi(t,y)\leq 0$, then
\begin{equation*}
y\leq \max\{y_1,y(0)\}\quad \textrm{on } (0,T).
\end{equation*}
\begin{proof}
%
See \cite[Lemma 3.3]{ChiyoDuzgunFrassuVigliaoro-2024}.
\end{proof}
\end{lemma}
\begin{lemma}\label{lemmap}
Let $n\in \N$, $m_2,m_3\in \R$, $\gamma,\alpha,\beta>0$, and  the lower bound of the relation \eqref{condgamma} be fixed. Then, for all $q>1$ there exists $\Bar{p}>1$ such that for all $p>\Bar{p}$ and
\begin{equation*}
    \theta(p)\coloneqq\frac{\frac{p+\gamma-1}{\gamma}-\frac{p+\gamma-1}{p\gamma}}{\frac{p+\gamma-1}{\gamma}-\frac 1\gamma+\frac 1n}, \quad \quad \sigma(p)\coloneqq\frac{p\gamma}{p+\gamma-1}, \quad \quad  \Bar{\theta}(p)\coloneqq \frac{\frac{p+\gamma-1}{\gamma}-\frac{p+\gamma-1}{\gamma(p+m_2+\alpha-1)}}{\frac{p+\gamma-1}{\gamma}-\frac 1\gamma+\frac 1n}, \quad \quad \Bar{\sigma}(p)\coloneqq\frac{\gamma(p+m_2+\alpha-1)}{p+\gamma-1},
\end{equation*}
\begin{equation*}
    \hat{\theta}(p)\coloneqq\frac{\frac{p+\gamma-1}{\gamma}-\frac{p+\gamma-1}{\gamma(p+m_3+\beta-1)}}{\frac{p+\gamma-1}{\gamma}-\frac 1\gamma+\frac 1n}, \quad \quad \hat{\sigma}(p)\coloneqq \frac{\gamma(p+m_3+\beta-1)}{p+\gamma-1}, \quad \quad \Tilde{\theta}(p) \coloneqq \frac{\frac{p+\gamma-1}{\gamma}-\frac{p+\gamma-1}{\gamma\beta}}{\frac{p+\gamma-1}{\gamma}-\frac 1\gamma+\frac 1n}, \quad \quad 
\end{equation*}
\begin{equation*}
    \underline{\theta}(p)\coloneqq \frac{\frac{p+\gamma-1}{\gamma}-\frac{p+\gamma-1}{\gamma q}}{\frac{p+\gamma-1}{\gamma}-\frac 1\gamma+\frac 1n}, \quad \quad \underline{\sigma}(p)\coloneqq \frac{\gamma(p+q)}{p+\gamma-1},
\end{equation*}
these inequalities hold:
\begin{table}[H]
\begin{subequations}
    \begin{subtable}[h]{0.19\textwidth}
        \centering
        \begin{equation}\label{theta}
        0<\theta<1,
        \end{equation}
    \end{subtable}
    \hfill
    \begin{subtable}[h]{0.19\textwidth}
        \centering
        \begin{equation}\label{sigmatheta}
        0<\frac{\sigma\theta}{\gamma}<1.
        \end{equation}
     \end{subtable}
    \hfill
      \begin{subtable}[h]{0.19\textwidth}
        \centering
       \begin{equation}\label{thetabar}
        0<\Bar{\theta}<1,
        \end{equation}
    \end{subtable}
    \hfill
    \begin{subtable}[h]{0.19\textwidth}
        \centering
        \begin{equation}\label{sigmathetabar}
        0<\frac{\Bar{\sigma}\Bar{\theta}}{\gamma}<1.
        \end{equation}
     \end{subtable}
     \hfill
     \begin{subtable}[h]{0.19\textwidth}
        \centering
       \begin{equation}\label{thetahat}
        0<\hat{\theta}<1,
        \end{equation}
    \end{subtable}
    \\
     \begin{subtable}[h]{0.19\textwidth}
        \centering
        \begin{equation}\label{sigmathetahat}
        0<\frac{\hat{\sigma}\hat{\theta}}{\gamma}<1.
        \end{equation}
     \end{subtable}
     \hfill
          \begin{subtable}[h]{0.19\textwidth}
        \centering
       \begin{equation}\label{thetatilde}
        0<\Tilde{\theta}<1,
        \end{equation}
    \end{subtable}
    \hfill
    \begin{subtable}[h]{0.19\textwidth}
        \centering
        \begin{equation}\label{sigmathetatilde}
        0<\frac{\hat{\sigma}\Tilde{\theta}}{\gamma}<1.
        \end{equation}
     \end{subtable}
     \hfill
          \begin{subtable}[h]{0.19\textwidth}
        \centering
       \begin{equation}\label{thetaunder}
        0<\underline{\theta}<1,
        \end{equation}
    \end{subtable}
    \hfill 
    \begin{subtable}[h]{0.19\textwidth}
        \centering
        \begin{equation}\label{sigmathetaunder}
    0<\frac{\underline{\sigma}\,\underline{\theta}}{\gamma}<1.
        \end{equation}
     \end{subtable}
\end{subequations}
\end{table}
\noindent In the specific \eqref{thetatilde} and \eqref{sigmathetatilde} are satisfied under the stronger assumption $\beta>1$.
\begin{proof}
By computing $\theta'(p)$, we observe that  $\theta'(p)>0$ for $p$ sufficiently large, so that $\theta(p)$ is definitely increasing. Moreover, since  $\lim_{p\to+\infty}{\theta(p)}=1$,  relation \eqref{theta} remains valid for some large $p$. Similarly we can obtain  \eqref{sigmatheta} (for which we use $\gamma>1$), \eqref{thetabar}, \eqref{thetahat}, whereas for \eqref{sigmathetabar} and \eqref{sigmathetahat} we exploit, respectively, $\gamma>\frac{n}{n+1}(m_2+\alpha)$ and $\gamma>\frac{n}{n+1}(m_3+\beta)$.  As to \eqref{thetaunder} and \eqref{sigmathetaunder}, we observe that
\begin{equation*}
\lim_{p\to+\infty}{\underline{\theta}(p)}=\lim_{p\to+\infty}{\frac{\underline{\sigma}(p)\underline{\theta}(p)}{\gamma}}=1-\frac 1q \in (0,1).
\end{equation*}
On the other hand, for $\beta>1$
\begin{equation*}
    \lim_{p\to+\infty}{\Tilde{\theta}(p)}=\lim_{p\to+\infty}{\frac{\hat{\sigma}(p)\Tilde{\theta}(p)}{\gamma}}=1-\frac 1\beta \in (0,1),
\end{equation*}
so that \eqref{thetatilde} and \eqref{sigmathetatilde} are fulfilled.

From the above reasoning, we can pick some $\bar{p}>1$ such that for all $p>\bar{p}$ the claim is given. 
\end{proof}
\end{lemma}
\section{Local well-posedness. Proof of Theorem \ref{theoremExistence}}
Let us establish existence and uniqueness of local classical solutions to problems \eqref{localproblem} and \eqref{nonlocalproblem}. In the specific, we will provide details for the local model (in its fully parabolic and parabolic-elliptic versions), being similar the reasoning concerning the
nonlocal one. It is worthwhile mentioning that oppositely to \cite[Theorem 1.1]{IshidaLankeitVigliloro-Gradient}, where the linear version of our problems is studied,  we cannot use estimates for Neumann heat kernels, and some more insights have to be faced. 
\subsubsection*{Proof of Theorem \texorpdfstring{\ref{theoremExistence}}{teo}}
Let us start with the \textit{existence} issue.  For any $R>0$,
let us consider for some $0<T\leq 1$, to be precised later on,  the closed convex
subset
$$S_{T}=\{0\le u \in C^{1,\frac{\delta}{2}}(\bar{\Omega}\times [0,T]):  \|u(\cdot,t)-u_0\|_{L^\infty(\Omega)}\le R,\, \; \text{for\ all}\
t\in[0,T]\}.$$
Once an element $\tilde{u}$ of $S_T$ is picked, from the properties of $f_1$ and $f_2$, one has that $f_1(\tilde{u}),f_1(\tilde{u})\in C^{\delta,\frac{\delta}{2}}(\bar{\Omega}\times [0,T])$; in the specific, the associated solution $v$ and $w$ to problems
 \begin{equation}\label{2.2}
		\begin{cases}
	\tau v_t-	\Delta v+ v=f_1(\tilde{u})&\text{in}\ \Omega\times(0,T),\\
		v_\nu=0&\text{on}\ \partial\Omega\times(0,T),\\
		\tau v(x,0)=\tau v_0(x),
	\end{cases}
\end{equation}
and
\begin{equation}\label{2.3}
	\begin{cases}
	\tau w_t-	\Delta w+ w=f_2(\tilde{u})&\text{in}\ \Omega\times(0,T),\\
		w_\nu=0&\text{on}\ \partial\Omega\times(0,T),\\
		\tau w(x,0)=\tau w_0(x),
	\end{cases}
\end{equation}
are, thanks to $\tau v_0,\tau w_0\in C_\nu^{2+\delta}(\bar{\Omega})$ and  classical elliptic and parabolic regularity results (\cite[Theorem 9.33]{BrezisBook},  \cite[Theorem 6.31]{GilbarTrudinger} for $\tau=0$, and \cite[Corollary 5.1.22]{LunardiBook}, \cite[Theorem IV.5.3]{LSUBookInequality} for $\tau=1$), such that 
\begin{equation}
v,w\in C^{2+\delta, \tau+\frac{\delta}{2}}(\bar\Omega \times [0,T])
\quad \textrm{and precisely}\quad \sup_{t\in [0,T]} \lVert v(\cdot, t)\rVert_{C^{2+\delta}(\bar\Omega)}+ \sup_{t\in [0,T]} \lVert w(\cdot, t)\rVert_{C^{2+\delta}(\bar\Omega)}\leq H,
\end{equation}
where $H=H(f_1,f_2,\tilde{u},\tau v_0,\tau w_0)>0$ is estimated by $\tilde{u}$, so essentially, $H=H(R)$.
Now, given these gained properties of $\nabla v$ and $\nabla w$, considered the smoothness of $\zeta \longmapsto (1+\zeta)^\vartheta$ for all $\vartheta \in \R$ and $\zeta\geq 0$, let us define 
\begin{equation}\label{PositionForLunardi}
\begin{cases}
\varphi_1: \Omega \times (0,T)\rightarrow \R^+, \textrm{ with }\varphi_1(x,t):= (\tilde{u}(x,t)+1)^{m_1-1},\\
\varphi_2: \Omega \times (0,T)\rightarrow \R^n, \textrm{ with }\varphi_2(x,t):= -\chi (\tilde{u}(x,t)+1)^{m_2-1}\nabla v(x,t)+  \xi (\tilde{u}(x,t)+1)^{m_3-1}\nabla w(x,t),\\
\varphi_3: \Omega \times (0,T)\times \R^{1+n}\rightarrow \R, \textrm{ with }\varphi(x,t,u, \nabla u):= \lambda u^\rho(x,t)-\mu u^k(x,t)-c|\nabla u(x,t)|^\gamma,
\end{cases}
\end{equation}
and also this system: 
\begin{equation}\label{2.1}
	\begin{cases}
		u_t=\nabla\cdot\left((\tilde{u}+1)^{m_1-1}\nabla u-\chi u(\tilde{u}+1)^{m_2-1}\nabla v
		+\xi u(\tilde{u}+1)^{m_3-1}\nabla w\right) +\lambda u^\rho -\mu u^k - c\abs{\nabla u}^\gamma
		&\text{in}\ \Omega\times(0,T),\\
		u_\nu=0&\text{on}\ \partial\Omega\times(0,T),\\
		u(x,0)=u_0(x)&x\in\overline{\Omega}.
	\end{cases}
\end{equation}
In the specific, with the positions in \eqref{PositionForLunardi}, problem \eqref{2.1} reads as
\begin{equation*}
	\begin{cases}
		u_t=\mathcal{A}u +\varphi(x,t,u,\nabla u)=\varphi_1 \Delta u +u \nabla \cdot \varphi_2+(\nabla \varphi_1+\varphi_2)\cdot \nabla u + \varphi_3
		&\text{in}\ \Omega\times(0,T),\\
		\mathcal{B}_1u=u_\nu=0&\text{on}\ \partial\Omega\times(0,T),\\
		u(x,0)=u_0(x)&x\in\overline{\Omega},
	\end{cases}
\end{equation*}
and it is seen that $\varphi$ complies with \cite[(7.3.5)]{LunardiBook}, exactly using the notation therein. As a consequence, by recalling that $u_0\in C^{2+\delta}_\nu(\bar\Omega)$,  \cite[Theorem 7.3.3]{LunardiBook} ensures the existence of some $0<T< 1$, as above, such that problem \eqref{2.1} has a unique solution 
\begin{equation*}
u\in C^{2+\delta,1+\frac{\delta}{2}}(\bar\Omega\times [0,T]).
\end{equation*} 
In particular,  this produces some positive constant $K=K(\tilde{u},\nabla v,\nabla w,u_0)$, and henceforth $K=K(R)$, with the property that
$$|u(x,t)-u_0(x)|\le Kt^{1+\tfrac{\delta}{2}} \quad \textrm{for all} \quad(x,t)\in\Omega\times(0,T), \quad \textrm{or}\quad \max_{t\in[0,T]}\,\|u(\cdot,t)-u_0\|_{L^\infty(\Omega)}\le KT^{1+\frac{\delta}{2}}. $$
In this way, we deduce 
that 
$$ \textrm{for} \quad T\leq \left(\frac{R}{K}\right)^\frac{2}{2+\delta}, \quad\|u(\cdot,t)-u_0\|_{L^\infty(\Omega)}\le R\;\quad
\textrm{for all} \quad t\in[0,T]. $$
Moreover, since $\underline{u}\equiv 0$ is a subsolution of the first equation in \eqref{2.1}, the parabolic comparison principle warrants the nonnegativity of $u$ on $\Omega \times (0,T)$. So, for values of $T$ as above, the map  $\Phi (\tilde{u})=u$ where $u$ solves problem \eqref{2.1}, is such that $\Phi(S_T) \subset S_T$ and $\Phi$ is compact, because the
 \cite[Ascoli--Arzel\`a Theorem 4.25]{BrezisBook} implies that the natural embedding
of $C^{2+\delta,\tau+\frac{\delta}{2}}(\bar{\Omega}\times [0,T]))$ into $C^{1,\frac{\delta}{2}}(\bar{\Omega}\times [0,T]))$ is a compact linear operator. Let $u$ be the fixed point of $\Phi$ asserted by the Schauder's fixed point theorem: first, the elliptic and parabolic maximum principles and $f_1(u),f_2(u)\geq 0$ in problems \eqref{2.2} and \eqref{2.3} also imply $v,w\geq 0$ in $\Omega \times (0,T)$; secondly,  as seen, $u,v,w$ have the claimed regularity. 

As to \textit{uniqueness}, let $(u_1,v_1,w_1)$ and $(u_2,v_2,w_2)$ be two different nonnegative classical solutions of problem \eqref{localproblem} in $\Omega\times (0,T)$ with the same initial data $u_1(\cdot,0)=u_2(\cdot,0)=u_0(x)$ and $\tau v_1(\cdot,0)=\tau v_2(\cdot,0)=\tau v_0(x)$ and $\tau w_1(\cdot,0)=\tau w_2(\cdot,0)=\tau w_0(x)$. By confining our attention to the case $\tau=0$ (being the case $\tau=1$ similar, as we will indicate later), we thus have from 
\eqref{2.3}--\eqref{2.1} and $i=1,2$  these six problems:
\begin{equation}\label{3.2_Bis}
	\begin{cases}
		-\Delta v_i+ v_i=f_1(u_i)&\text{in}\ \Omega\times(0,T),\\
		(v_i)_\nu=0&\text{on}\ \partial\Omega\times(0,T),
	\end{cases}
\end{equation}
\begin{equation}\label{3.3_Bis}
	\begin{cases}
		-\Delta w_i+ w_i=f_2(u_i)&\text{in}\ \Omega\times(0,T),\\
		(w_i)_\nu=0&\text{on}\ \partial\Omega\times(0,T),
	\end{cases}
\end{equation}
and
\begin{equation}\label{3.1_Bis}
	\begin{cases}
 \!\begin{aligned}[t]
  &(u_i)_t=\nabla\cdot ((u_i+1)^{m_1-1}\nabla u_i
		-\chi u_i(u_i+1)^{m_2-1}\nabla v_i
		+\xi u_i(u_i+1)^{m_3-1}\nabla w_i)\\
  &\quad \quad \quad \quad \quad  +\lambda u_i^\rho-\mu u_i^k-c|\nabla u_i|^\gamma  
\end{aligned} 
 & \text{in}\ \Omega\times(0,T), 
\\
		(u_i)_\nu=0&\text{on}\ \partial\Omega\times(0,T),
  \\
		u_i(x,0)=u_0(x)&x\in\overline{\Omega}.
	\end{cases}
\end{equation}
 Now we set 
\begin{align*}
	s_1=\min\left\{\|u_1\|_{L^\infty(\Omega\times(0,T))},
	\|u_2\|_{L^\infty(\Omega\times(0,T))}\right\},\quad
	s_2=\max\left\{\|u_1\|_{L^\infty(\Omega\times(0,T))},
	\|u_2\|_{L^\infty(\Omega\times(0,T))}\right\},
\end{align*}
and introduce for $s\ge0$ the function
$$\theta(s)=\begin{dcases}\log(s+1)&\text{if}\ m_1=0,\\
	\frac{(s+1)^{m_1}}{m_1}&\text{if}\ m_1\neq 0,\end{dcases}
\quad \textrm{with }\; m_1\in\R.$$
Noting that $\theta^\prime(s)>0$ for all $s\ge0$, we define, in view of  forthcoming applications of the Mean Value Theorem, the positive $T$-dependent constants 
\begin{equation}\label{ConstantsC1-2-3}
	\begin{dcases}
	C_1=\max_{[s_1,s_2]}|f_1^\prime|,\\ C_2=\max_{[s_1,s_2]}|f_2^\prime|,
\end{dcases}
\;\textrm{ and }\;
	C_3=\begin{dcases}
		\max_{[s_1,s_2]}\theta^\prime&\text{if}\ (u_1-u_2)\Delta(u_1-u_2)\geq 0,\\
		\min_{[s_1,s_2]} \theta^\prime&\text{if}\ (u_1-u_2)\Delta(u_1-u_2)\leq 0.
	\end{dcases}
\end{equation}
%
Under such circumstances, by considering problems \eqref{3.2_Bis}, we have that $(v_1-v_2)$ solves, for some $\bar{u}\in(u_1,u_2)$, the equation $-\Delta(v_1-v_2)+(v_1-v_2)=f_1(u_1)-f_1(u_2)=f_1'(\bar{u})(u_1-u_2)$, with Neumann boundary conditions. In light of this, through Young's inequality and taking into account \eqref{ConstantsC1-2-3}, testing against $(v_1-v_2)$, the mentioned equation provides the following estimate
\begin{align*}
	\int_\Omega|\nabla(v_1-v_2)|^2&+\int_\Omega(v_1-v_2)^2=\int_\Omega(f_1(u_1)-f_1(u_2))
	(v_1-v_2)\leq C_1\int_\Omega |(u_1-u_2)(v_1-v_2)|\\ &
	\le\int_\Omega(v_1-v_2)^2
	+\const{Giu-2}\int_\Omega(u_1-u_2)^2 \quad \textrm{on }\, t\in(0,T),
\end{align*}
inferring
\begin{equation}\label{3.4}
	\int_\Omega|\nabla(v_1-v_2)|^2\le
	\const{Giu-2}\int_\Omega(u_1-u_2)^2\quad \textrm{for all }\, t\in(0,T).
\end{equation}
Similarly, by considering problems \eqref{3.2_Bis} we arrive at
\begin{equation}\label{3.4BISS}
\int_\Omega|\nabla(w_1-w_2)|^2\le
\const{Giu-3}\int_\Omega(u_1-u_2)^2\quad \textrm{for all }\, t\in(0,T).
\end{equation}
We also get from problems \eqref{3.1_Bis}
\begin{equation}\label{MainForUniqueness}
\begin{split}
	\frac{d\mathcal{F}}{dt}=\frac{1}{2}\frac{d}{dt}\int_\Omega(u_1-u_2)^2&=\int_\Omega(u_1-u_2)(u_1-u_2)_t\\
	&=-\int_\Omega\nabla(u_1-u_2) \cdot 
	\left((u_1+1)^{m_1-1}\nabla u_1-(u_2+1)^{m_1-1}\nabla u_2\right)\\
	&\quad +\chi\int_\Omega\nabla(u_1-u_2)\cdot 
	\left(u_1(u_1+1)^{m_2-1}\nabla v_1-u_2(u_2+1)^{m_2-1}\nabla v_2\right)\\
	&\quad -\xi\int_\Omega\nabla(u_1-u_2)\cdot 
	\left(u_1(u_1+1)^{m_3-1}\nabla w_1-u_2(u_2+1)^{m_3-1}\nabla w_2\right)\\
	& \quad +\int_\Omega (u_1-u_2)(\lambda u_1^\rho-\lambda u_2^\rho-\mu u_1^k+\mu u_2^k-c|\nabla u_1|^\gamma+c|\nabla u_2|^\gamma) \quad \textrm{ for all } \,t\in (0,T),
	\end{split}
\end{equation}
where the first, second, third, and fourth integrals at the right-hand side are denoted, for simplicity, by $I_1$, $I_2$, $I_3$ and $I_4$, respectively. We now estimate these separate terms as follows (recall \eqref{ConstantsC1-2-3}, again, and the definition of $\theta$):
\begin{equation}
		\label{3.6}
\begin{split}
	I_1&=-\int_\Omega\nabla(u_1-u_2)\cdot \nabla\left(\theta(u_1)-\theta(u_2)\right)=\int_\Omega\Delta(u_1-u_2)(\theta(u_1)-\theta(u_2))
	=\int_\Omega\Delta(u_1-u_2)\theta^\prime(\overline{u})(u_1-u_2)\\
	&\le C_3\int_\Omega(u_1-u_2)\Delta(u_1-u_2)=-C_3\int_\Omega|\nabla(u_1-u_2)|^2 \quad \textrm{for all } t\in(0,T).
\end{split}
\end{equation}
On the other hand, once more  the smoothness of $\zeta \longmapsto (1+\zeta)^
\vartheta$ for all $\vartheta \in \R$ and $\zeta\geq 0$, makes that for some $C_4$ and $C_5$, and for 
$j\in\{m_2,m_3\}$
$$
\begin{cases}
	|u_1(u_1+1)^{j-1}-u_2(u_2+1)^{j-1}|\le C_4|u_1-u_2|
	&\text{in}\; \Omega\times(0,T),\\
	u_2(u_2+1)^{j-1}\le C_4 &\text{in} \;\Omega\times(0,T);
	\end{cases}
 $$
 furthermore
 $$
 \begin{cases}
	|\nabla v_1|\le C_5\ \text{and}\ |\nabla v_2|\le C_5&\text{in}\ \Omega\times(0,T),\\
	|\nabla w_1|\le C_5\ \text{and}\ |\nabla w_2|\le C_5
	&\text{in}\ \Omega\times(0,T).\end{cases}$$
For similar reasons, there exists $C_6$ such that for all $\rho \geq 1$
\begin{equation}
|u_1^\rho-u_2^\rho| \leq C_6|u_1-u_2| \quad \textrm{ and } \quad ||\nabla u_1|^\rho-|\nabla u_2|^\rho|\leq C_6||\nabla u_1|-|\nabla u_2|| \quad \text{ in}\ \Omega\times(0,T).
\end{equation} 
Henceforth, starting from the identity
\begin{align*}
	I_2=\chi\int_\Omega\nabla(u_1-u_2)&\cdot \left((u_1(u_1+1)^{m_2-1}-u_2(u_2+1)^{m_2-1})
	\nabla v_1+u_2(u_2+1)^{m_2-1}\nabla(v_1-v_2)\right) \quad \textrm{ on } (0,T),
\end{align*}
and using Hölder's inequality and relation  \eqref{disab} we obtain
\begin{align*}
	I_2^2&\le\chi^2\int_\Omega|(u_1(u_1+1)^{m_2-1}-u_2(u_2+1)^{m_2-1})\nabla v_1+u_2(u_2+1)^{m_2-1}\nabla(v_1-v_2)|^{2}
	\int_\Omega|\nabla(u_1-u_2)|^2\\
	&\le2\chi^2\left(\int_\Omega|(u_1(u_1+1)^{m_2-1}-u_2(u_2+1)^{m_2-1})\nabla v_1|^2
	+\int_\Omega|u_2(u_2+1)^{m_2-1}\nabla(v_1-v_2)
	|^2\right)\int_\Omega|\nabla(u_1-u_2)|^2
	\\
	&\le2\chi^2C_4^2\left(C_5^2\int_\Omega(u_1-u_2)^2
	+\int_\Omega|\nabla(v_1-v_2)|^2\right)\int_\Omega|\nabla(u_1-u_2)|^2 \quad \textrm{ for all } t\in (0,T).
\end{align*}
Taking the square root, we get by Young's inequality
\begin{equation}
	\label{3.7}
\begin{split}
	I_2&\le\chi C_4\sqrt{2}\!\left(C_5^2\int_\Omega(u_1\!-\!u_2)^2
	\!+\!\int_\Omega|\nabla(v_1\!-\!v_2)|^2\right)^{\tfrac{1}{2}}
	\!\!\left(\int_\Omega|\nabla(u_1\!-\!u_2)|^2\right)^{\tfrac{1}{2}}\\
	&\le\frac{C_3}{3}\int_\Omega|\nabla(u_1-u_2)|^2
	+\const{Giu-0}\int_\Omega(u_1-u_2)^2
	+\const{Giu-0}\int_\Omega|\nabla(v_1-v_2)|^2 \quad \textrm{ on }  (0,T).
\end{split}
\end{equation}
Similarly, from 
\begin{align*}
	I_3=\xi\int_\Omega&  \nabla(u_1-u_2)\cdot \left(\left(u_2(u_2+1)^{m_3-1}-u_1(u_1+1)^{m_3-1}\right)
	\nabla w_2+u_1(u_1+1)^{m_3-1}\nabla(w_2-w_1)\right) \quad \textrm{for every } t\in (0,T),
\end{align*}
 we obtain
 \begin{equation}\label{3.8}
 	\begin{split}
 		I_3&\le\xi C_5\sqrt{2}\!\left(C_5^2\int_\Omega(u_1\!-\!u_2)
 		\!+\!\int_\Omega|\nabla(w_1\!-\!w_2)|^2\right)^{\tfrac{1}{2}}
 		\!\left(\int_\Omega|\nabla(u_1\!-\!u_2)|^2\right)^{\tfrac{1}{2}}
 		\\ &
 		\leq \frac{C_3}{3}\int_\Omega|\nabla(u_1-u_2)|^2+\const{Giu-23}
 		\int_\Omega(u_1-u_2)^2 +\const{Giu-23}\int_\Omega |\nabla (w_1-w_2)|^2 \quad \textrm{ on }\, (0,T).
 	\end{split}
 \end{equation}
%
As to the addendum
\begin{equation}
I_4:=\int_\Omega (u_1-u_2)(\lambda u_1^\rho-\lambda u_2^\rho-\mu u_1^k+\mu u_2^k-c|\nabla u_1|^\gamma+c|\nabla u_2|^\gamma) \quad \textrm{on }  (0,T),
\end{equation}
analogous computations based on the Mean Value Theorem lead to 
\begin{equation}
I_4\leq C_6 \int_\Omega (u_1-u_2)^2+C_6\int_\Omega |u_1-u_2||\nabla u_1||-|\nabla u_2|| \quad \textrm{on }  (0,T),
\end{equation}
and in turn Young's inequality, using that $(|\nabla u_1|-|\nabla u_2|)^2\leq |\nabla (u_1-u_2)|^2$, provides 
\begin{equation}\label{FinaleUnicita}
I_4\leq \const{GV}\int_\Omega (u_1-u_2)^2+\frac{C_3}{3}\int_\Omega |\nabla (u_1-u_2)|^2 \quad \textrm{for all }  t\in (0,T).
\end{equation}
Now a rearrangement of the derived bounds \eqref{3.4}--\eqref{FinaleUnicita} yields
a suitable positive constant $\tilde{C}$ such that 
$$\mathcal{F}'(t):=\frac{1}{2}\frac{d}{dt}\int_\Omega (u_1-u_2)^2  \le\tilde{C}\int_\Omega (u_1-u_2)^2=\tilde{C}\mathcal{F}(t),\qquad t\in(0,T).$$ We then exploit the nonnegativity of $\mathcal{F}(t)$, the initial
condition $\mathcal{F}(0)=0$ (based on $u_1(\cdot,0)=u_2(\cdot,0)$), and Lemma \ref{LemmaODI-Comparison} to prove that $\mathcal{F}(t)\equiv0$. As a result, $u_1-u_2=0$ on $\Omega \times (0,T)$ and  consequently, recalling again problems \eqref{3.2_Bis} and \eqref{3.3_Bis}, we manifestly get $v_1-v_2=0$ and $w_1-w_2=0$, as well on $\Omega \times (0,T)$; the proof is concluded. 
\qed
\begin{remark}
    The proof of Theorem \ref{theoremExistence} when model \eqref{nonlocalproblem} is considered, is essentially justified by the same steps than those used for problem \eqref{localproblem} when $\tau=0$. Specifically, as to the uniqueness of the solution $v$ and $w$ for the Poisson equations under homogeneous Neumann boundary conditions, this is ensured from $\int_\Omega v(x,t)dx=\int_\Omega w(x,t)dx=0$. (Recall Remark \ref{RemarkDeviation}.) For $\tau=1$, the idea is based on the evolutive behavior of $\mathcal{F}(t)=\frac{1}{2}\int_\Omega \left( (u_1-u_2)^2+(v_1-v_2)^2+(w_1-w_2)^2\right).$ 
\end{remark}
\section{Dichotomy and boundedness criterion. The requirement on the growth of the gradient term}
From Theorem \ref{theoremExistence}, we have these consequences: 
\begin{corollary}\label{existence}  Under the hypotheses of Theorem \ref{theoremExistence}, there exists $\TM \in (0,\infty]$ such that the solution  $(u,v,w)$ to problems \eqref{localproblem} and \eqref{nonlocalproblem} obeys this dichotomy criterion:
 \begin{equation}\label{dictomyCriteC2+del} 
 \text{either }\, \TM=\infty,\, \text{ or} \quad \TM<\infty \quad \text{and then} \quad \limsup_{t \to \TM} \left(\|u(\cdot,t)\|_{C^{2+\delta}(\bar\Omega)}+\|v(\cdot,t)\|_{C^{2+\delta}(\bar\Omega)}+\|w(\cdot,t)\|_{C^{2+\delta}(\bar\Omega)}\right)=\infty.
 \end{equation}
 Finally, if additionally $k>\rho$ the mass is uniformly bounded in time, in the sense that 
\begin{equation}\label{boundednessMass}
  u\in L^\infty((0,\TM);L^1(\Omega)), \textrm{ and precisely } \int_{\Omega} u \leq M:=\max\left\{\int_\Omega u_0(x)dx, \left(\frac{\lambda}{\mu} |\Omega|^{k-\rho}\right)^{\frac{1}{k-\rho}}\right\} \quad \text{for all } t \in [0,\TM).   
\end{equation} 
\begin{proof} 
 Taking $T$ as initial time and $u(\cdot, T)$ as the initial condition, Theorem \ref{theoremExistence} would provide a solution $(\hat{u},\hat{v},\hat{w})$ defined on $\bar{\Omega}\times [T,\hat{T}]$, for some $\hat{T}>0$, which by uniqueness would be the prolongation of $(u,v,w)$ (exactly, from $\bar{\Omega}\times [0,T]$ to  $\bar{\Omega}\times [0,\hat{T}])$. This procedure can be repeated up to construct a maximal interval time $[0,\TM)$ of existence, in the sense that either $\TM=\infty$, or if $\TM<\infty$ no solution belonging to $C^{2+\delta,\tau+\frac{\delta}{2}}(\bar{\Omega}\times [0,\TM])$ may exist and henceforth relation \eqref{dictomyCriteC2+del} has to be fulfilled. 
 
 An integration of the first equation in problem \eqref{localproblem} and an application of the divergence theorem give, thanks to the Neumann boundary conditions,
\begin{equation}\label{mprimo}
\begin{split}
    \frac{d}{dt}\int_\Omega u =:y'\leq\lambda\into u^\rho-\mu\into u^k \quad \text{in }(0,\TM).
\end{split}
\end{equation}
Now, thanks to the H\"{o}lder inequality, we have for $k>\rho\geq 1$
\begin{equation}\label{uk}
    -\into u^k \leq -\abs*{\Omega}^{\frac{\rho-k}{\rho}}\tonda*{\into u^\rho}^\frac{k}{\rho} \quad \textrm{and} \quad -\into u^\rho \leq -\abs*{\Omega}^{1-\rho}\tonda*{\into u}^\rho\quad \text{in }(0,\TM).
\end{equation}
By combining expressions \eqref{mprimo} and \eqref{uk} we arrive for $\gamma(t)=\lambda \int_\Omega u^\rho(x,t)dx\geq 0$ on $(0,\TM)$ at this initial problem
\begin{equation*}
\begin{dcases}
    y'\leq \gamma(t)\left(1 - \frac{\mu}{\lambda} \abs*{\Omega}^{\rho-k} y^{k-\rho}\right) & \text{in }(0,\TM),\\
    y(0)=\into u_0(x)dx,
\end{dcases}
\end{equation*}
so concluding by invoking Lemma \ref{LemmaODI-Comparison} with $$T=\TM,\, \phi(t,y)=\gamma(t)\left(1 - \frac{\mu}{\lambda} \abs*{\Omega}^{\rho-k} y^{k-\rho}\right),\, y(0)=\into u_0(x)dx\, \textrm{ and }\, y_1=\left(\frac{\lambda}{\mu} |\Omega|^{k-\rho}\right)^{\frac{1}{k-\rho}}.$$
\end{proof}
\end{corollary}
From now on with $(u,v,w)$ we will refer to the unique local solutions to models \eqref{localproblem} and \eqref{nonlocalproblem} defined on $\Omega \times (0,\TM)$ and  provided by Corollary \ref{existence}. 
As seen in such a corollary, relation \eqref{dictomyCriteC2+del} holds and in this case the solution is said to blow-up at $\TM$ in $C^{2+\delta}(\Omega)$-norm; nevertheless, this does not necessarily imply the blow-up  in $L^\infty(\Omega)$-norm. Exactly in order to rephrase the extensibility criterion in \eqref{dictomyCriteC2+del} in terms of some uniform-in-time boundedness of $u$, we have to restrict the class of gradient nonlinearities of the type $|\nabla u|^\gamma$ for $\gamma \geq 1$ to those whose growth is at the most quadratic. We precisely have this 
\begin{lemma}[Extension and boundedness criteria] \label{ExtensionLemma}Let $\gamma \in [1,2]$. If $u\in L^\infty((0,\TM);L^\infty(\Omega))$, then we have $\TM=\infty$ and in particular $u\in L^\infty((0,\infty);L^\infty(\Omega))$.
\begin{proof}
By contradiction let $\TM$ be finite. Naturally, it is sufficient to show that boundedness of $u$ in $L^\infty(\Omega)$ for $t\in(0,\TM)$, entails finiteness of $\lVert u(\cdot,t) \rVert_{C^{2+\delta}(\bar{\Omega})}+\lVert v(\cdot,t) \rVert_{C^{2+\delta}(\bar{\Omega})}+\lVert w(\cdot,t) \rVert_{C^{2+\delta}(\bar{\Omega})}$ on $[0,\TM].$ In fact, from \eqref{dictomyCriteC2+del}, we would have an inconsistency. 

First, we have boundedness of $\lVert v(\cdot,t) \rVert_{C^{1+\delta}(\bar{\Omega})}$, and of $\lVert w(\cdot,t) \rVert_{C^{1+\delta}(\bar{\Omega})}$ on $[0,\TM]$, as consequence of our hypothesis $u\in L^\infty((0,\TM);L^\infty(\Omega))$, once \cite[Theorem~1.2]{lieberman_paper} (if $\tau=1$) or, e.g., from \cite[Theorem ~I.19.1]{friedman_pde} (if $\tau=0$) and the Sobolev embedding immersion $W^{2,p}(\Omega)\subset C^{1+\delta}(\bar{\Omega})$, with arbitrarily large $p$, are invoked. On the other hand, according to the nomenclature in \cite{lieberman_paper}, letting $X=(x,t)$, $A(X,u,\nabla u)=-(u+1)^{m_1-1}\nabla u+\chi u(u+1)^{m_2-1}\nabla v(x,t)-\xi u(u+1)^{m_3-1}\nabla w(x,t) \textrm{ and } B(X,u,\nabla u)=-\lambda u^\rho+\mu u^k+c|\nabla u|^\gamma,$ we conclude from \cite[Theorem ~1.2]{lieberman_paper} that $\lVert u(\cdot,t) \rVert_{C^{1+\delta}(\bar{\Omega})}$ is finite as well on $[0,\TM]$, by using that $u$ is a bounded solution to
\begin{equation*}
\begin{cases}
u_t+\nabla \cdot A(X,u,\nabla u)+B(X,u,\nabla u)=0 & \textrm{ in } \Omega \times (0,\TM),\\
A(X,u,\nabla u)\cdot \nu=0 & \textrm{ on } \partial \Omega \times (0,\TM),\\
u(x,0)=u_0(x)& x\in \bar{\Omega},
\end{cases}
\end{equation*}
exactly for $\gamma\in [1,2]$ (see in the specific \cite[($1.2c$),  Theorem ~1.2]{lieberman_paper}). From this gained regularity of $u$ and the relative bound,  bootstrap  arguments already developed in Corollary \ref{existence}, in conjunction with the regularity of the initial data $(u_0,\tau v_0,\tau w_0)$ for model \eqref{localproblem} give 
$$
\sup_{t\in [0,\TM]}\left(\lVert u(\cdot,t) \rVert_{C^{2+\delta}(\bar{\Omega})}+\lVert v(\cdot,t) \rVert_{C^{2+\delta}(\bar{\Omega})}+\lVert w(\cdot,t) \rVert_{C^{2+\delta}(\bar{\Omega})}\right)<\infty.
$$
\end{proof}
\end{lemma} 

The last result indicates how to ensure uniform-in-time boundedness of $u$ on $(0,\TM)$ exploiting some uniform-in-time estimate of $u$ on $(0,\TM)$ in some $L^p(\Omega)$-norm.
\begin{lemma}\label{extensioncriterion}
Let $f_1(u)$ and $f_2(u)$ comply with \eqref{disf} and $\gamma \in [1,2]$. If $u\in L^\infty((0,\TM);L^p(\Omega))$ for every $p>1$, then $u\in L^\infty((0,\infty);L^\infty(\Omega))$.
\begin{proof}
With the aid of Lemma \ref{ExtensionLemma}, it is sufficient to show  that  $u\in L^\infty((0,\TM);L^\infty(\Omega))$. First we observe that for any $\lambda,\mu,c>0$, and $k>\rho \geq 1$, $\gamma\geq 1$, the function $u \mapsto \lambda u^\rho-\mu u^{k}-c|\nabla u|^\gamma$, for $u\geq 0$, is bounded from above by a positive constant $L$;
henceforth, we can write $g(x,t):=\lambda u^\rho-\mu u^{k}-c|\nabla u|^\gamma\leq L$ in $\Omega \times (0,\TM)$. Subsequently, with a view to \cite[Appendix A]{TaoWinkParaPara}, and consistently to that symbology, it is seen that if $u$  classically solves in $\Omega \times (0,\TM)$ the first equation in \eqref{localproblem} (and \eqref{nonlocalproblem}), it also solves \cite[problem (A.1)]{TaoWinkParaPara} with 
\begin{equation*}
D(x,t,u)=(u+1)^{m_1-1},\quad f(x,t)=-\chi u(u+1)^{m_2-1}\nabla v+\xi  u(u+1)^{m_3-1}\nabla w,\quad g(x,t)=L. 
 \end{equation*}
In particular, this in conjunction with the Neumann boundary conditions on $v$ and $w$ imply that hypotheses (A.2)--(A.5) are met, whereas the second inclusion of (A.6) holds for any choice of $q_2>1$. Since we are assuming $p$ arbitrarily large, and  $u\in L^\infty((0,\TM);L^p(\Omega))$, conditions (A.7)--(A.10) are automatically complied. Let us now deal with the first inclusion of \cite[(A.6)]{TaoWinkParaPara}. If $\tau=0$ in problem \eqref{localproblem} (or directly considering model \eqref{nonlocalproblem}), from $u\in L^\infty((0,\TM);L^p(\Omega))$ we have $f_1(u)\in L^\infty((0,\TM);L^p(\Omega))$ so that elliptic regularity results, and embeddings, imply  $\nabla v\in L^\infty((0,\TM);L^\infty(\Omega))$. Since the same inclusion holds for $\nabla w$, we have that, for any $q_1>1$,
\begin{equation*}
    f=-\chi u(u+1)^{m_2-1}\nabla v +\xi u(u+1)^{m_3-1}\nabla w \in L^\infty((0,\TM);L^{q_1}(\Omega)).
\end{equation*}
If $\tau=1$, we can reason in the same way invoking \eqref{tau1extension}. So we have the claim exactly by virtue of \cite[Lemma A.1]{TaoWinkParaPara}.
\end{proof}
\end{lemma}
\section{Derivation of sufficient a priori estimates. Proof of Theorem \ref{theoremlocal}}
The following two lemmas will be used later on in the paper. (The appearing value $\bar{p}$ is taken from Lemma \ref{lemmap}.)
\begin{lemma} Under the hypotheses of Lemma \ref{lemmap}, for every $\varepsilon$ we have  that for all $p>\Bar{p}$
\begin{equation}\label{disp}
    \into(u+1)^p\leq \varepsilon \into\abs*{\nabla (u+1)^\frac{p+\gamma-1}{\gamma}}^\gamma + \const{ss} \quad \text{on }(0,\TM),
\end{equation}
\begin{equation}\label{disphifinal}
    -\into\abs*{\nabla (u+1)^\frac{p+\gamma-1}{\gamma}}^\gamma\leq \const{sd} - \const{sf}\tonda*{\into(u+1)^p}^\frac{\gamma}{\sigma\theta} \quad \text{on }(0,\TM),
\end{equation}
\begin{equation}\label{dis1}
    \into(u+1)^{p+m_2+\alpha-1}\leq \varepsilon \into\abs*{\nabla (u+1)^\frac{p+\gamma-1}{\gamma}}^\gamma + \const{sg} \quad \text{on }(0,\TM),
\end{equation}
\begin{equation}\label{dis2}
    \into(u+1)^{p+m_3+\beta-1}\leq \varepsilon\into\abs*{\nabla (u+1)^\frac{p+\gamma-1}{\gamma}}^\gamma + \const{sjj} \quad \text{on }(0,\TM),
\end{equation}
\begin{equation} \label{u+1w}
    \into (u+1)^{p+m_3-1}w \leq \varepsilon\into (u+1)^{p+m_3+\beta-1}+ \varepsilon\into \abs*{\nabla (u+1)^{\frac{p+\gamma-1}{\gamma}}}^\gamma + \const{sj1} \quad\text{on }(0,\TM), 
\end{equation}
\begin{equation} \label{intmol}
    \into (u+1)^{p+m_3-1}\into (u+1)^\beta \leq \varepsilon\into (u+1)^{p+m_3+\beta-1}+ \varepsilon\into \abs*{\nabla (u+1)^{\frac{p+\gamma-1}{\gamma}}}^\gamma + \const{sj2} \quad\text{on }(0,\TM).
\end{equation}
\begin{proof}
In order to prove relation \eqref{disp}, we use the Gagliardo--Nirenberg inequality (see \cite{Nirenber_GagNir_Ineque}) in conjunction with \eqref{theta} and \eqref{boundednessMass}, obtaining
\begin{equation*}
\begin{split}
    \into (u+1)^p=&\norm*{(u+1)^\frac{p+\gamma-1}{\gamma}}_{L^{\frac{p\gamma}{p+\gamma-1}}(\Omega)}^{\frac{p\gamma}{p+\gamma-1}}\\ 
    \leq& \const{d_1}\tonda*{\norm*{\nabla(u+1)^\frac{p+\gamma-1}{\gamma}}_{L^\gamma(\Omega)}^\theta \norm*{(u+1)^\frac{p+\gamma-1}{\gamma}}_{L^\frac{\gamma}{p+\gamma-1}(\Omega)}^{1-\theta}+\norm*{(u+1)^\frac{p+\gamma-1}{\gamma}}_{L^\frac{\gamma}{p+\gamma-1}(\Omega)}}^\sigma\\
    \leq& \const{d_3}\tonda*{\into\abs*{\nabla (u+1)^\frac{p+\gamma-1}{\gamma}}^\gamma}^\frac{\sigma\theta}{\gamma}+\const{d_4}\leq \varepsilon \into\abs*{\nabla (u+1)^\frac{p+\gamma-1}{\gamma}}^\gamma + \const{ss} \quad \text{on }(0,\TM),
\end{split}
\end{equation*}
where in the second-last step we employed relation \eqref{disab}, and in the last one Young's inequality, justified by relation \eqref{sigmatheta}.

As to bound \eqref{disphifinal} we refrain by using the Young inequality in the previous derivation: subsequently we have
\begin{equation*}
\begin{split}
    \into (u+1)^p=&\norm*{(u+1)^\frac{p+\gamma-1}{\gamma}}_{L^{\frac{p\gamma}{p+\gamma-1}}(\Omega)}^{\frac{p\gamma}{p+\gamma-1}}
    \leq \const{d_3}\tonda*{\into\abs*{\nabla (u+1)^\frac{p+\gamma-1}{\gamma}}^\gamma}^\frac{\sigma\theta}{\gamma}+\const{d_4}
    \quad \text{on }(0,\TM).
\end{split}
\end{equation*}
and we conclude again by relying on \eqref{disab}.

On the other hand, relation \eqref{dis1} is achieved again by invoking the Gagliardo--Nirenberg inequality, \eqref{thetabar}, the boundedness of the mass (i.e., \eqref{boundednessMass}), Young's inequality and \eqref{sigmathetabar}: specifically we have
\begin{equation*}
\begin{split}
    \into (u+1)^{p+m_2+\alpha-1}=&\norm*{(u+1)^\frac{p+\gamma-1}{\gamma}}_{L^{\frac{\gamma(p+m_2+\alpha-1)}{p+\gamma-1}}(\Omega)}^\frac{\gamma(p+m_2+\alpha-1)}{p+\gamma-1}\\
    \leq& \const{s2}\tonda*{\norm*{\nabla(u+1)^\frac{p+\gamma-1}{\gamma}}_{L^\gamma(\Omega)}^{\Bar{\theta}}  \norm*{(u+1)^\frac{p+\gamma-1}{\gamma}}_{L^\frac{\gamma}{p+\gamma-1}(\Omega)}^{1-\Bar{\theta}}+\norm*{(u+1)^\frac{p+\gamma-1}{\gamma}}_{L^\frac{\gamma}{p+\gamma-1}(\Omega)}}^{\Bar{\sigma}}\\
    \leq& \const{s3}\tonda*{\into\abs*{\nabla (u+1)^\frac{p+\gamma-1}{\gamma}}^\gamma}^\frac{\Bar{\sigma}\Bar{\theta}}{\gamma}+\const{s4} \leq \varepsilon \into\abs*{\nabla (u+1)^\frac{p+\gamma-1}{\gamma}}^\gamma + \const{sg}\quad \text{on }(0,\TM).
\end{split}
\end{equation*}
Moreover, in this circumstance relations  \eqref{thetahat} and \eqref{sigmathetahat} give
\begin{equation*}
\begin{split}
    \into (u+1)^{p+m_3+\beta-1}=&\norm*{(u+1)^\frac{p+\gamma-1}{\gamma}}_{L^{\frac{\gamma(p+m_3+\beta-1)}{p+\gamma-1}}(\Omega)}^{\frac{\gamma(p+m_3+\beta-1)}{p+\gamma-1}}\\ 
    \leq& \const{dd_1}\tonda*{\norm*{\nabla(u+1)^\frac{p+\gamma-1}{\gamma}}_{L^\gamma(\Omega)}^{\hat{\theta}} \norm*{(u+1)^\frac{p+\gamma-1}{\gamma}}_{L^\frac{\gamma}{p+\gamma-1}(\Omega)}^{1-\hat{\theta}}+\norm*{(u+1)^\frac{p+\gamma-1}{\gamma}}_{L^\frac{\gamma}{p+\gamma-1}(\Omega)}}^{\hat{\sigma}}\\
    \leq& \const{dd_3}\tonda*{\into\abs*{\nabla (u+1)^\frac{p+\gamma-1}{\gamma}}^\gamma}^\frac{\hat{\sigma}\hat{\theta}}{\gamma}+\const{dd_4}  \leq \varepsilon \into\abs*{\nabla (u+1)^\frac{p+\gamma-1}{\gamma}}^\gamma + \const{sjj} \quad \text{for all } t\in (0,\TM). 
\end{split}
\end{equation*}
We also note, by using Young's inequality that 
\begin{equation*}
    \into (u+1)^{p+m_3-1}w \leq \varepsilon\into(u+1)^{p+m_3+\beta-1} +\const{zzzO}\into w^{\frac{p+m_3+\beta-1}{\beta}} \quad \text{on }(0,\TM),
\end{equation*}
so that from \eqref{ellipticreg} and \eqref{disf}
    \begin{equation}\label{PD}
    \begin{split}
        \into (u+1)^{p+m_3-1}w & \leq \varepsilon\into(u+1)^{p+m_3+\beta-1}+\varepsilon\into f_2(u)^{\frac{p+m_3+\beta-1}{\beta}}+\const{oopp}\tonda*{\into f_2(u)}^{\frac{p+m_3+\beta-1}{\beta}}\\
        &\leq \varepsilon\into(u+1)^{p+m_3+\beta-1}+\const{asf2}\tonda*{\into (u+1)^\beta}^{\frac{p+m_3+\beta-1}{\beta}} \quad \text{on }(0,\TM),
    \end{split}
    \end{equation}
    where $\const{asf2}=\const{asf2}(C_\varepsilon).$     Now, in order to estimate $\const{asf2}\tonda*{\into (u+1)^\beta}^{\frac{p+m_3+\beta-1}{\beta}}$, we obtain that if $\beta\leq1$, relation \eqref{boundednessMass}  implies
 \begin{equation*}
       \const{asf2}\tonda*{\into (u+1)^\beta}^{\frac{p+m_3+\beta-1}{\beta}}
       \leq  \const{sj1}\quad \text{on }(0,\TM).
    \end{equation*}
    Conversely, for $\beta>1$, by properly using the Gagliardo--Nirenberg and Young inequalities, relations \eqref{thetatilde} and \eqref{sigmathetatilde}, and \eqref{boundednessMass}, we can write
    \begin{equation*}
    \begin{split}
       \const{asf2}\tonda*{\into (u+1)^\beta}^{\frac{p+m_3+\beta-1}{\beta}}&=\const{asf2}\norm*{(u+1)^{\frac{p+\gamma-1}{\gamma}}}_{L^{\frac{\beta\gamma}{p+\gamma-1}}(\Omega)}^{\frac{\gamma(p+m_3+\beta-1)}{p+\gamma-1}} \\
       &\leq \const{asf4}\tonda*{\norm*{\nabla(u+1)^{\frac{p+\gamma-1}{\gamma}}}^{\Tilde{\theta}}_{L^{\gamma}(\Omega)}\norm*{(u+1)^{\frac{p+\gamma-1}{\gamma}}}^{1-\Tilde{\theta}}_{L^{\frac{\gamma}{p+\gamma-1}}(\Omega)}+\norm*{(u+1)^{\frac{p+\gamma-1}{\gamma}}}_{L^{\frac{\gamma}{p+\gamma-1}}(\Omega)}   }^{\hat{\sigma}}\\
       &\leq \const{asf5}\tonda*{\into\abs*{\nabla (u+1)^{\frac{p+\gamma-1}{\gamma}}}^\gamma}^{\frac{\hat{\sigma}\Tilde{\theta}}{\gamma}}+\const{asf6}\\
       &\leq \varepsilon\into \abs*{\nabla (u+1)^{\frac{p+\gamma-1}{\gamma}}}^\gamma +\const{sj1} \quad \text{on }(0,\TM).
    \end{split}
    \end{equation*}
    Coming back to bound \eqref{PD}, the two derivations above established show that for $\beta>0$  estimate \eqref{u+1w} is valid.
    
    We similarly operate to derive \eqref{intmol}: if $\beta\leq1$, it is directly obtained by using  the boundedness of the mass (i.e., \eqref{boundednessMass}) and then Young's inequality  so to arrive  at
    \begin{equation*}
        \into (u+1)^{p+m_3-1}\into (u+1)^\beta\leq \const{ad1}\into (u+1)^{p+m_3-1}\leq \varepsilon\into (u+1)^{p+m_3+\beta-1} + \const{sj2}\quad \text{on }(0,\TM).
    \end{equation*}
    If $\beta>1$, we fix $q>\max{\{\beta,m_3+\beta-1\}}$ and we get through the Hölder inequality and \eqref{LemmaEsponenti} that
    \begin{equation}\label{PM}
    \begin{split}
    \into (u+1)^{p+m_3-1}\into (u+1)^\beta &\leq \const{zs3}\quadra*{\into(u+1)^{p+m_3+\beta-1} }^{\frac{p+m_3-1}{p+m_3+\beta-1}} \quadra*{\tonda*{\into (u+1)^q}^{\frac pq +1}}^{\frac{\beta}{p+q}}   \\
    & \leq \varepsilon \into(u+1)^{p+m_3+\beta-1} + \tonda*{\into (u+1)^q}^{\frac pq +1} +\const{as6}\quad \textrm{on} \;(0,\TM).
    \end{split}
    \end{equation}
 On the other hand, by virtue of $u\in L^\infty((0,\TM);L^1(\Omega))$ and relations \eqref{thetaunder} and  \eqref{sigmathetaunder} we can write 
    \begin{equation}\label{PDDP}
    \begin{split}
    \tonda*{\into (u+1)^q}^{\frac pq +1}=&\norm*{(u+1)^\frac{p+\gamma-1}{\gamma}}_{L^{\frac{\gamma q}{p+\gamma-1}}(\Omega)}^\frac{\gamma(p+q)}{p+\gamma-1}\\
    \leq& \const{ss2}\tonda*{\norm*{\nabla(u+1)^\frac{p+\gamma-1}{\gamma}}_{L^\gamma(\Omega)}^{\underline{\theta}}  \norm*{(u+1)^\frac{p+\gamma-1}{\gamma}}_{L^\frac{\gamma}{p+\gamma-1}(\Omega)}^{1-\underline{\theta}}+\norm*{(u+1)^\frac{p+\gamma-1}{\gamma}}_{L^\frac{\gamma}{p+\gamma-1}(\Omega)}}^{\underline{\sigma}}\\
    \leq& \const{sd3}\tonda*{\into\abs*{\nabla (u+1)^\frac{p+\gamma-1}{\gamma}}^\gamma}^\frac{\underline{\sigma}\underline{\theta}}{\gamma}+\const{ss4} \leq \varepsilon \into\abs*{\nabla (u+1)^\frac{p+\gamma-1}{\gamma}}^\gamma + \const{ss42}\quad \text{on }(0,\TM). 
    \end{split}
    \end{equation}
By inserting \eqref{PDDP} into \eqref{PM} also for $\beta>1$ bound \eqref{intmol} is achieved. 
\end{proof}
\end{lemma}
In the previous result, we did not use that $(u,v,w)$ is a solution of the models we are considering. The forthcoming lemma, indeed, is valid for solutions to both problems \eqref{localproblem} and \eqref{nonlocalproblem}. 
\begin{lemma}\label{validoperentrambi}Let the same hypotheses of Lemma \ref{lemmap} be given. 
Then 
\begin{equation}\label{disprincipale}
    \varphi'(t)\leq-\const{b}\into F_{m_2}(u)\Delta v +\const{c}\into F_{m_3}(u)\Delta w -\const{e}\into\abs*{\nabla(u+1)^{\frac{p+\gamma-1}{\gamma}}}^\gamma + \const{cc} \quad \text{on }(0,\TM),
\end{equation}
where $\varphi(t)=\frac{1}{p}\into (u+1)^p$ and, for $j\in\{m_2,m_3\}$
\begin{equation}\label{funcor}
    F_j(u) := \int_0^u \hat{u} \tonda*{\hat{u}+1}^{p+j-3} d \hat{u}.
\end{equation}
\begin{proof}
Let us differentiate the functional $\varphi(t)=\frac 1p \into (u+1)^p$. Using the first equation of \eqref{localproblem} and the divergence theorem, we have for every $t\in(0,\TM)$
\begin{equation*}
\begin{split}
    \varphi'(t)=&\into (u+1)^{p-1}u_t \\ =&\into (u+1)^{p-1}\nabla\cdot ((u+1)^{m_1-1}\nabla u)-\chi\into (u+1)^{p-1}\nabla\cdot (u(u+1)^{m_2-1}\nabla v)\\
    &+\xi\into(u+1)^{p-1}\nabla\cdot (u(u+1)^{m_3-1}\nabla w) +\lambda\into (u+1)^{p-1}u^\rho-\mu\into (u+1)^{p-1}u^k-c\into (u+1)^{p-1}\abs{\nabla u}^\gamma\\
    \leq&-(p-1)\into(u+1)^{p+m_1-3}\abs{\nabla u}^2+\chi (p-1)\into u(u+1)^{p+m_2-3}\nabla u\cdot\nabla v -\xi (p-1)\into u(u+1)^{p+m_3-3}\nabla u\cdot \nabla w \\ &+\lambda\into (u+1)^{p+\rho-1}-\mu\into (u+1)^{p-1}u^k-c\left(\frac{p-1+\gamma}{\gamma}\right)^{-\gamma}\into \abs*{\nabla(u+1)^{\frac{p+\gamma-1}{\gamma}}}^\gamma.
\end{split}
\end{equation*}
By considering the definition of $F_j(u)$ in \eqref{funcor}, we have $\nabla F_j(u)=u(u+1)^{p+j-3}\nabla u$, so that  again the divergence theorem yields on $(0,\TM)$
\begin{equation*}
    \varphi'(t)\leq-\const{b}\into F_{m_2}(u)\Delta v +\const{c}\into F_{m_3}(u)\Delta w +\lambda\into (u+1)^{p+\rho-1}-\mu\into (u+1)^{p-1}u^k-\const{e}\into\abs*{\nabla(u+1)^{\frac{p+\gamma-1}{\gamma}}}^\gamma.
\end{equation*}
Inverting relation \eqref{disab}  we also have 
\begin{equation*}
    -\mu\into (u+1)^{p-1}u^k \leq -\frac{\mu}{2^{k-1}}\into (u+1)^{p+k-1}+\mu\into (u+1)^{p-1} \quad \textrm{for all } t<\TM,
\end{equation*}
and consequently 
\begin{equation}\label{IneqQuasiUlti}
\begin{split}
    \varphi'(t)\leq&-\const{b}\into F_{m_2}(u)\Delta v +\const{c}\into F_{m_3}(u)\Delta w +\lambda\into (u+1)^{p+\rho-1}-\const{gh8}\into (u+1)^{p+k-1}\\&+\mu\into (u+1)^{p-1}-\const{e}\into\abs*{\nabla(u+1)^{\frac{p+\gamma-1}{\gamma}}}^\gamma \quad\text{on }(0,\TM).
\end{split}
\end{equation}
Moreover, since $k>\rho\geq 1$, we can invoke Young's inequality to observe that
\begin{equation}\label{disprho}
   \lambda \into (u+1)^{p+\rho-1} \leq \frac{\const{gh8}}{2} \into (u+1)^{p+k-1}+\const{df} \quad\text{on }(0,\TM),
\end{equation}
and 
\begin{equation}\label{disp1}
   \mu  \into (u+1)^{p-1} \leq \frac{\const{gh8}}{2}\into (u+1)^{p+k-1}+\const{df1} \quad\text{on }(0,\TM).
\end{equation}
Finally, by using bounds \eqref{disprho} and \eqref{disp1}  into relation \eqref{IneqQuasiUlti} we obtain what claimed.
\end{proof}
\end{lemma}
\subsection{The local model: analysis of problem \texorpdfstring{\eqref{localproblem}}{L}} In this section we will give the proof of Theorem \ref{theoremlocal}, 
by deriving the crucial uniform-in-time boundedness of the solution $u$ in some $L^p(\Omega)$-norm. We will separate the cases $\tau=0$ and $\tau=1$. 

Supported by assumptions \eqref{disf}, we observe that from relation \eqref{funcor} we have
\begin{equation}\label{disF}
    \frac{u^{p+j-1}}{p+j-1}
 \leq F_j(u)\leq \frac{1}{p+j-1}\left[(u+1)^{p+j-1}-1\right].
\end{equation}
We will make use of this property.
\subsubsection{The parabolic-elliptic case: problem \texorpdfstring{$(\mathcal{L}_0)$}{L0}}
Let us start with the case $\tau=0$.
\begin{lemma}\label{lemma2local}
Let the hypotheses of Lemma \ref{lemmap} be given. 
Then for all $p>\bar{p}$ we have that  $u\in L^\infty((0,\TM);L^p(\Omega))$. 
\begin{proof}
By taking into account \eqref{disprincipale} and \eqref{disF}, by using the second and third equations of \eqref{localproblem} we obtain on $(0,\TM)$
\begin{equation*}
\varphi'(t)\leq\const{zx1}\into(u+1)^{p+m_2+\alpha-1}+\const{zx2}\into(u+1)^{p+m_3-1}w -\const{zx3}\into u^{p+m_3-1}f_2(u)-\const{e}\into\abs*{\nabla(u+1)^{\frac{p+\gamma-1}{\gamma}}}^\gamma + \const{cc}.
\end{equation*}
Successively, by exploiting \eqref{disab} and \eqref{disf}, we get on $(0,\TM)$
\begin{equation*}
\begin{split}
\varphi'(t)\leq&\const{zx1}\into(u+1)^{p+m_2+\alpha-1}+\const{zx2}\into (u+1)^{p+m_3-1}w-\const{zxz1}\into(u+1)^{p+m_3+\beta-1}+\const{0T}\int_\Omega (u+1)^\beta\\& \quad -\const{e}\into \abs*{\nabla(u+1)^{\frac{p+\gamma-1}{\gamma}}}^\gamma+\const{cc},
\end{split}
\end{equation*}
and by manipulating $\const{0T}\int_\Omega (u+1)^\beta$ through the Young inequality we have on $(0,\TM)$
\begin{equation*}
\begin{split}
\varphi'(t)\leq&\const{zx1}\into(u+1)^{p+m_2+\alpha-1}+\const{zx2}\into (u+1)^{p+m_3-1}w-\const{asdf1}\into(u+1)^{p+m_3+\beta-1}-\const{e}\into \abs*{\nabla(u+1)^{\frac{p+\gamma-1}{\gamma}}}^\gamma+\const{a_67}.
\end{split}
\end{equation*}
In this position,  by invoking \eqref{dis1} and \eqref{u+1w} for proper small  $\varepsilon$ we have
\begin{equation*}
    \varphi'(t)\leq -\const{a_5s}\into \abs*{\nabla(u+1)^{\frac{p+\gamma-1}{\gamma}}}^\gamma+\const{a_6s} \quad \text{on $(0,\TM)$},
\end{equation*}
so that an  application of \eqref{disphifinal} yields
\begin{equation}\label{problemainiziale}
\begin{dcases}
    \varphi'(t)\leq \const{k1}-\const{l1}\varphi(t)^{\frac{\gamma}{\sigma\theta}} \quad \textrm{for all } t\in(0,\TM),\\
    \varphi(0)=\frac 1p \into (u_0+1)^p.
\end{dcases}
\end{equation}
Lemma \ref{LemmaODI-Comparison} infers $\varphi(t)\leq \max\left\{\varphi(0),\left(\frac{\const{k1}}{\const{l1}}\right)^{\frac{\sigma\theta}{\gamma}}\right\}$ for all $0<t<\TM$, i.e., $u\in L^\infty((0,\TM);L^p(\Omega)).$
\end{proof}
\end{lemma}
\subsubsection{The parabolic-parabolic case: problem \texorpdfstring{$(\mathcal{L}_1)$}{L1}} Let us continue the anlysis of model \eqref{localproblem} when $\tau=1$. We will make use of the parabolic regularity result given in Lemma \ref{regularity}.
\begin{lemma}\label{lemma3local} Let the hypotheses of Lemma \ref{lemmap} be given. Moreover let $f_1(u)$ and $f_2(u)$ comply with \eqref{disf}.
 Then for all $p>\bar{p}$ we have that  $u\in L^\infty((0,\TM);L^p(\Omega))$. 
\begin{proof}
In the claim of Lemma \ref{validoperentrambi},  let us estimate the integrals involving $\Delta v$ and $\Delta w$  by invoking relation \eqref{disF} and Young's inequality. We have  for all $t\in(0,\TM)$ 
\begin{equation}\label{second}
    \const{b}\into F_{m_2}(u)\abs{\Delta v}\leq \const{b_2}\into (u+1)^{p+m_2-1}\abs{\Delta v}\leq \into (u+1)^{p+m_2+\alpha-1} + \const{b_4}\into \abs{\Delta v}^\frac{p+m_2+\alpha-1}{\alpha}
\end{equation}
and
\begin{equation}\label{third}
    \const{c}\into F_{m_3}(u)\abs{\Delta w}\leq \const{c_2}\into (u+1)^{p+m_3-1}\abs{\Delta w}\leq \into (u+1)^{p+m_3+\beta-1} + \const{c_4}\into \abs{\Delta w}^\frac{p+m_3+\beta-1}{\beta}.
\end{equation}
By plugging into bound \eqref{disprincipale} estimates \eqref{second} and \eqref{third}, we obtain
\begin{equation}\label{dislemma1local}
\begin{split}
    \varphi'(t)\leq & \into(u+1)^{p+m_2+\alpha-1} + \const{b_4}\into \abs{\Delta v}^\frac{p+m_2+\alpha-1}{\alpha} +\into(u+1)^{p+m_3+\beta-1}\\ +&\const{c_4}\into \abs{\Delta w}^\frac{p+m_3+\beta-1}{\beta}-\const{e}\into \abs*{\nabla(u+1)^{\frac{p+\gamma-1}{\gamma}}}^\gamma+ \const{cc} \quad \text{for all $t\in(0,\TM)$.}
    \end{split}
\end{equation}
 With relation \eqref{tau1} in our hands, by considering $z=v$, $g=f_1(u)$ (and using in addition \eqref{disf}) and $q=\frac{p+m_2+\alpha-1}{\alpha}$, we derive
\begin{equation}\label{regparv}
    \int_0^t e^s \into \abs*{\Delta v}^{\frac{p+m_2+\alpha-1}{\alpha}}\,ds \leq \const{fg} +  \const{fg}\int_0^t e^s \into(u+1)^{p+m_2+\alpha-1}\,ds \quad \text{for all } t\in(0,\TM),
\end{equation}
and, similarly for the equation of $w$,
\begin{equation}\label{regparw}
    \int_0^t e^s \into \abs*{\Delta w}^{\frac{p+m_3+\beta-1}{\beta}}\,ds \leq \const{dfs1} + \const{dfs1}\int_0^t e^s \into(u+1)^{p+m_3+\beta-1}\,ds \quad \text{for all } t\in(0,\TM),
\end{equation}
with $\const{fg}=\const{fg}(C_P)$, $\const{dfs1}=\const{dfs1}(C_P)$. In turn, relation \eqref{dislemma1local} can be also be rephrased as
\begin{equation*}
\begin{split}
    \varphi'(t)+ \varphi(t)\leq & \into(u+1)^{p+m_2+\alpha-1} + \const{b_4}\into \abs{\Delta v}^\frac{p+m_2+\alpha-1}{\alpha} +\into(u+1)^{p+m_3+\beta-1}\\ &+\const{c_4}\into \abs{\Delta w}^\frac{p+m_3+\beta-1}{\beta} +\frac 1p\into (u+1)^p -\const{e}\into \abs*{\nabla(u+1)^{\frac{p+\gamma-1}{\gamma}}}^\gamma+\const{cc} \quad \text{for all $t\in(0,\TM)$,}
\end{split}
\end{equation*}
and by solving this differential inequality we have on $(0,\TM)$,
\begin{equation*}
\begin{split}
    e^t \varphi(t) \leq& \varphi(0) +\int_0^t e^s \biggl(\into(u+1)^{p+m_2+\alpha-1} + \const{b_4}\into \abs{\Delta v}^\frac{p+m_2+\alpha-1}{\alpha}\\ & +\into(u+1)^{p+m_3+\beta-1}+\const{c_4}\into \abs{\Delta w}^\frac{p+m_3+\beta-1}{\beta} +\frac 1p \into (u+1)^p -\const{e}\into \abs*{\nabla(u+1)^{\frac{p+\gamma-1}{\gamma}}}^\gamma+\const{cc} \biggr) \,ds.
\end{split}
\end{equation*}
Now relations \eqref{regparv} and \eqref{regparw} entail on $(0,\TM)$
\begin{equation*}
\begin{split}
    \frac{1}{p}e^t \into u^p \leq  \frac{1}{p}e^t \into (u+1)^p =e^t\varphi(t)\leq& \const{gh1} +\int_0^t e^s \biggl( \const{a_1s}\into(u+1)^{p+m_2+\alpha-1}\\ & +\const{a_3s}\into(u+1)^{p+m_3+\beta-1}+\frac 1p\into (u+1)^p -\const{e}\into \abs*{\nabla(u+1)^{\frac{p+\gamma-1}{\gamma}}}^\gamma+ \const{cc} \biggr) \,ds.
\end{split}
\end{equation*}
We now use \eqref{disp}, \eqref{dis1} and \eqref{dis2},  with small enough $\varepsilon$, in a such way that the gradient term in the above inequality may  control the remaining three integrals so to arrive at
\begin{equation*}
    e^t \into u^p\leq  \const{1gh1} + \const{1a_6e}(e^t -1) \quad \textrm{on }  (0,\TM),
\textrm{ or equivalently }
    \into u^p\leq \const{sdff2}  \quad \textrm{for all } t \in (0,\TM).
\end{equation*}
\end{proof}
\end{lemma}
\subsection{The nonlocal model: analysis of problem \texorpdfstring{\eqref{nonlocalproblem}}{NL}}
Let us now turn our attention to problem \eqref{nonlocalproblem}, for which the goal is  establishing for $u$ the same inclusion obtained in Lemma \ref{lemma2local} as well.
\begin{lemma}\label{lemmanonlocal}Let the hypotheses of Lemma \ref{lemmap} be fixed. Moreover let $f_1(u)$ and $f_2(u)$ comply with \eqref{disf}. Then for all $p>\bar{p}$ we have that $u\in L^\infty((0,\TM);L^p(\Omega))$.
\begin{proof}
Let us exploit again relation \eqref{disprincipale}; we use the second and third equations of model \eqref{nonlocalproblem} so that, by neglecting nonpositive terms, we have on $(0,\TM)$ and taking in mind restrictions in  \eqref{disf}
\begin{equation*}
\begin{split}
    \varphi'(t)\leq &-\const{b}\into F_{m_2}(u)\Delta v +\const{c}\into F_{m_3}(u)\Delta w -\const{e}\into\abs*{\nabla(u+1)^{\frac{p+\gamma-1}{\gamma}}}^\gamma+\const{cc}\\
    \leq & \const{b}\into F_{m_2}(u)f_1(u) +\const{cw}\into F_{m_3}(u)\into f_2(u)-\const{c}\into F_{m_3}(u)f_2(u)-\const{e}\into\abs*{\nabla(u+1)^{\frac{p+\gamma-1}{\gamma}}}^\gamma+\const{cc}\\
    \leq & \const{bww}\into (u+1)^{p+m_2+\alpha-1} +\const{cww}\into (u+1)^{p+m_3-1}\into (u+1)^\beta-\const{cw2}\into u^{p+m_3-1}(u+1)^\beta-\const{e}\into\abs*{\nabla(u+1)^{\frac{p+\gamma-1}{\gamma}}}^\gamma+\const{cc}.
\end{split}
\end{equation*}
Similarly to what done previously, by using relation \eqref{disab} and the Young inequality on the third term of the right-hand side of the now obtained estimate, we have
\begin{equation*}
    \varphi'(t)\leq \const{bww}\into (u+1)^{p+m_2+\alpha-1} +\const{cww}\into (u+1)^{p+m_3-1}\into (u+1)^\beta-\const{cw22}\into (u+1)^{p+m_3+\beta-1}-\const{e}\into\abs*{\nabla(u+1)^{\frac{p+\gamma-1}{\gamma}}}^\gamma+\const{ssaa}.
\end{equation*}
By recalling relations \eqref{dis1} and \eqref{intmol}, we again have for proper $\varepsilon$ and up to constants problem \eqref{problemainiziale}, and we conclude.
\end{proof}
\end{lemma}

\subsubsection*{Proof of Theorem \ref{theoremlocal}}
Depending to which case one is referring, we use Lemma \ref{lemma2local} (for problem \eqref{localproblem} with $\tau=0$) or Lemma \ref{lemma3local} (for problem \eqref{localproblem} with $\tau=1$) or Lemma \ref{lemmanonlocal} (for problem \eqref{nonlocalproblem})  and, possibly enlarging $\Bar{p}$, Lemma \ref{extensioncriterion} gives the claim.

\section{Numerical simulations}

\subsection{Numerical approximation and tools}
The numerical approximation of chemotaxis models is far to be  straightforward, especially when dealing with blow-up situations, which are very challenging to reproduce in the discrete framework. In fact, any positive and mass-bounded discrete solution, in the sense that the $\norm{\cdot}_{L^1(\Omega)}$ is bounded, will never show a proper blow-up phenomenon as it is also bounded in norm $\norm{\cdot}_{L^\infty(\Omega)}$ due to the equivalence of the norms in a discrete, finite-dimensional, space. Moreover, these chemotactic collapse situations lead to the appearance of very steep gradients which may produce severe nonphysical spurious oscillations as shown in \cite{chertock_second-order_2008,gutierrez-santacreu_analysis_2021}. In this work these difficulties can even be strengthen due to the damping gradient term introduced in systems \eqref{localproblem} and \eqref{nonlocalproblem} that may tend to infinity when a chemotactic collapse occurs.

In order to overcome these inconveniences and provide an accurate and meaningful numerical approximation for the solution of the related mechanism, we are inspired by the approach in \cite{acosta-soba_KS_2022}, where the classical Keller--Segel model is rewritten as a gradient flow. In particular, regarding the spatial approximation, we use an upwind discontinuous Galerkin (DG) method \cite{di2011mathematical} with piecewise constant polynomials defined on a triangular/tetrahedral mesh of the ($2$-dimensional/$3$-dimensional) domain $\Omega$ with size $h$. This is combined with a proper implicit-explicit approximation in time, defined on an equispaced partition of the temporal domain with time step $\Delta t$. By means of this method, we can obtain a positive approximating solution of the problem which prevents spurious oscillations from appearing. Moreover, the resulting scheme is linear, and it decouples the equations for the chemical signals, $v$ and $w$, which are computed first in each time step, from the equation for the cell density, $u$, computed afterwards in each time step. 

As to the employed software, we have used the Python interface of the open source library FEniCSx \cite{ScroggsEtal2022}.
\subsection{The general setting and a blow-up criterion. Simulating model \texorpdfstring{$(\mathcal{L}_0)$}{L0}}
\begin{enumerate}
\item We confine our numerical studies to the parabolic-elliptic-elliptic version of problem \eqref{localproblem}, since we expect the results to be similar when $\tau=1$ and for model \eqref{nonlocalproblem}. (Some tests in similar attraction-repulsion contexts for the case $\tau=1$ can be found, for instance, in \cite{FrassuGalvanViglialoro}.)
\item We define $f_1(u)=u^\alpha$, $f_2(u)=u^\beta$ and we assume that all the parameters are set to $1$, but $k=1.1$, unless otherwise specified. For the sake of clarity, we also indicate in the figures any different value of the parameters with respect to those already fixed.
\item As to the detection of blow-up phenomena, since an accurate discrete solution idealizing such scenario will exhibit a mass accumulation in small regions of the domain, leading to an exponential growth in a certain time interval, the following \textit{blow-up criterion} appears consistent (see  \cite{shen_unconditionally_2020,badia2023bound,acosta-soba_KS_2022} for details):
\begin{criterion}[Blow-up criterion]\label{CriterionBU}
Blow-up occurs when the norm $\|\cdot\|_{L^\infty(\Omega)}$  of the approximated solution stabilizes over time. An estimate of the  blow-up time $\TM$ is the approximated value of the instant of time around which the stabilization exhibits.
\end{criterion}
\end{enumerate}

\subsubsection{Numerical simulations for the linear and only attraction model}

\begin{figure}
    \centering
    \begin{subfigure}[b]{0.32\textwidth}
    \centering
        \includegraphics[width=\textwidth]{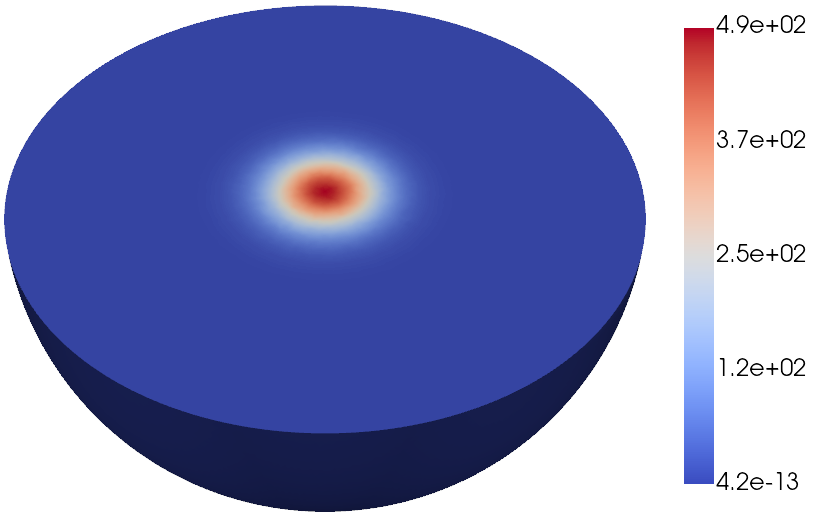}
        \caption{$t=0$.}
        \label{fig:test_1_p1c_u_a}
    \end{subfigure}
    \hfill
    \begin{subfigure}[b]{0.32\textwidth}
    \centering
        \includegraphics[width=\textwidth]{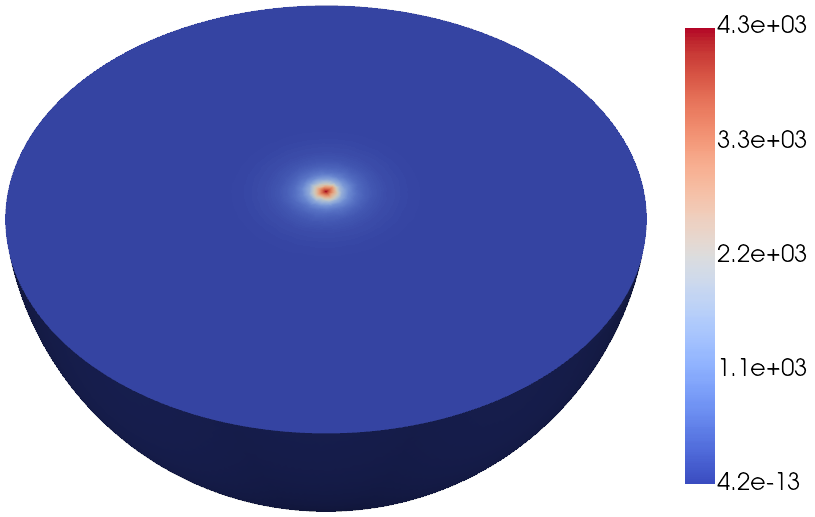}
        \caption{$t=5\cdot 10^{-4}$.}
        \label{fig:test_1_p1c_u_b}
    \end{subfigure}
    \hfill
    \begin{subfigure}[b]{0.32\textwidth}
    \centering
        \includegraphics[width=\textwidth]{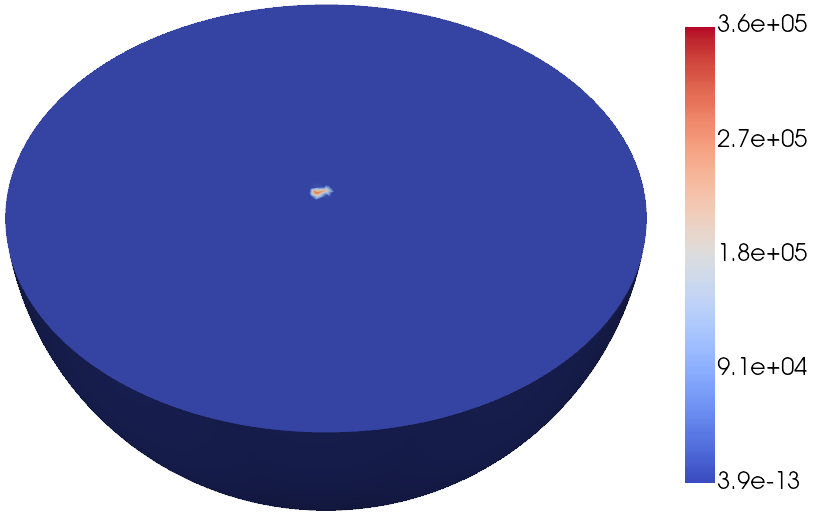}
        \caption{$t=3\cdot 10^{-3}$.}
        \label{fig:test_1_p1c_u_c}
    \end{subfigure}      
    \caption{Approximation of the solution $u$ at different time steps ($c=0$, $\chi=5$, $\xi=0$).}
    \label{fig:test_1_p1c_u}
\end{figure}

\begin{figure}
    \centering
    \begin{subfigure}[b]{0.49\textwidth}
    \centering
        \includegraphics[width=\textwidth]{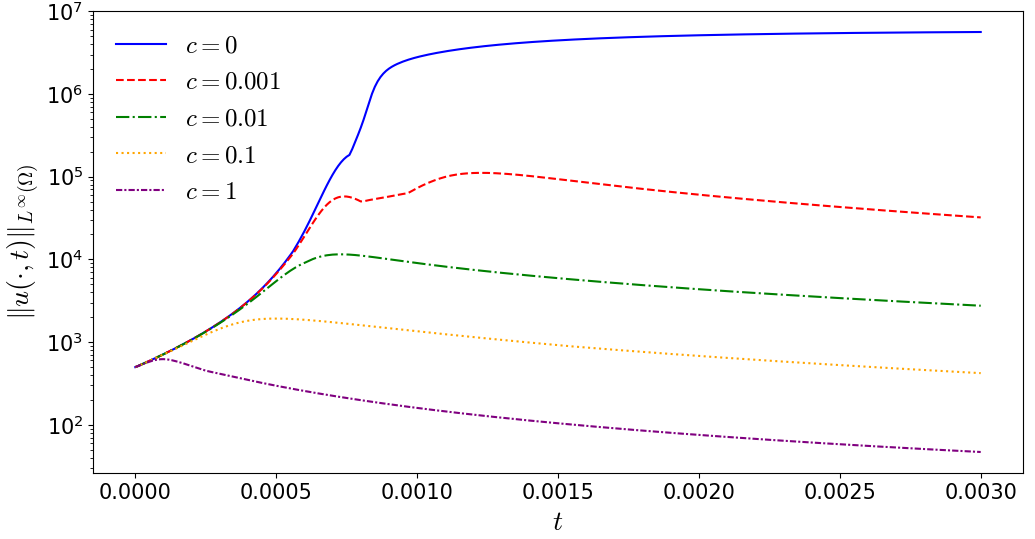}
        \caption{$\gamma=1.75$}
        \label{fig:test_1_max-u_a}
    \end{subfigure}
    \hfill
    \begin{subfigure}[b]{0.49\textwidth}
    \centering
        \includegraphics[width=\textwidth]{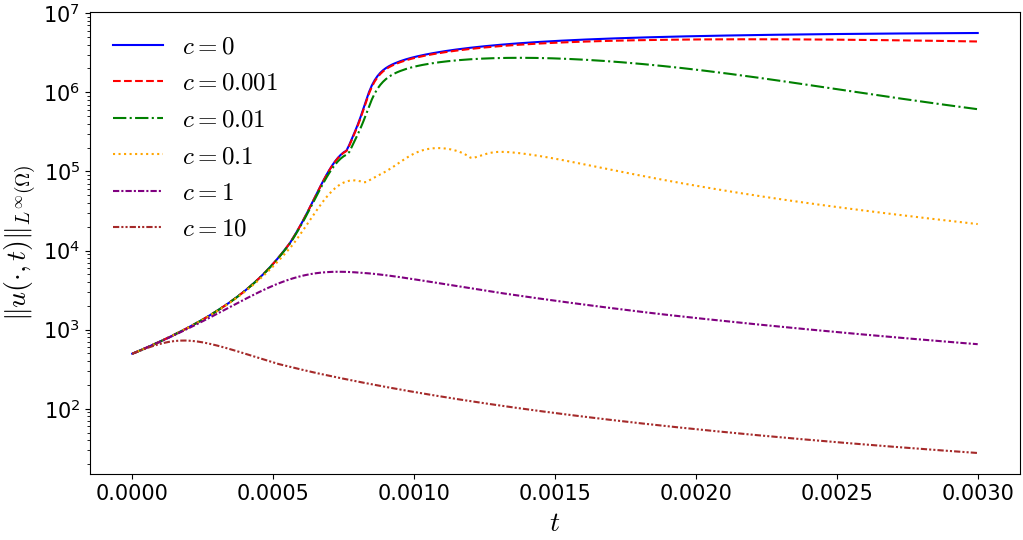}
        \caption{$\gamma=1.4$.}
        \label{fig:test_1_max-u_b}
    \end{subfigure}
    \\
    \begin{subfigure}[b]{0.49\textwidth}
    \centering
        \includegraphics[width=\textwidth]{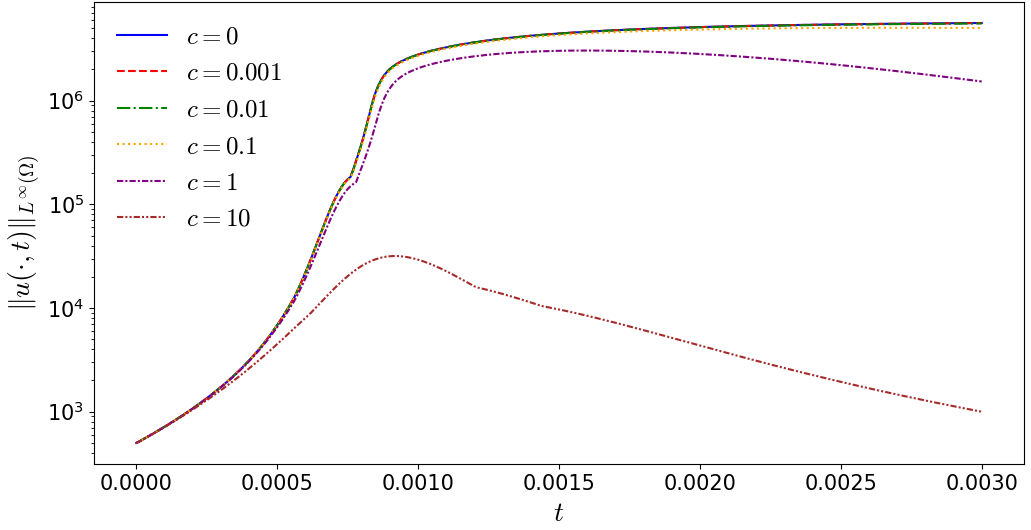}
        \caption{$\gamma=1.1$.}
        \label{fig:test_1_max-u_c}
    \end{subfigure}      
    \caption{Maximum of the approximation of $u$ over time for different values of $c$ and $\gamma$ ($\chi=5$, $\xi=0$).}
    \label{fig:test_1_max-u}
\end{figure}

In this first example, we consider a 3D version ($n=3$) of model \eqref{localproblem} defined in the spatial domain $\Omega=\{(x,y,z)\in\R^3\colon x^2+y^2+z^2< 1\}$ with the following initial condition
$$
u_0(x,y,z)\coloneqq500e^{-35(x^2+y^2+z^2)},
$$
plotted in Figure~\ref{fig:test_1_p1c_u_a}. Moreover, we take $\chi=5$ and $\xi=0$ (no repulsion) and we use a mesh of size $h\approx4.4\cdot 10^{-2}$ and a time step $\Delta t=10^{-5}$. In Figure~\ref{fig:test_1_p1c_u} the approximation obtained with $c=0$ is plotted at different time steps. We can observe that a blow-up phenomenon seems to appear (indeed, it is observable in Figures \ref{fig:test_1_p1c_u_b} and  \ref{fig:test_1_p1c_u_c} how  the maximum of the solution increases in the center of the sphere, eventually achieving the value $10^5$), according with the analytical results given in \cite{Winkler_ZAMP-FiniteTimeLowDimension} -- notice that $k=1.1<7/6\approx 1.167$ so the condition  in \cite[(1.4), Theorem 1.1]{Winkler_ZAMP-FiniteTimeLowDimension} is satisfied. This blow-up phenomenon (interpretable in the sense of Criterion \ref{CriterionBU} in the next figures) can be prevented using appropriate damping gradient terms, as established in Theorem~\ref{theoremlocal} (generalization of \cite[Theorem 1.2]{IshidaLankeitVigliloro-Gradient}). As a matter of fact, we show in Figure~\ref{fig:test_1_max-u} the evolution of the maximum of the approximations for different choices of $c$ and $\gamma$. Notice that the maximum shown in Figure~\ref{fig:test_1_max-u} ($c=0$) achieves even higher values, up to the order $10^6$, than those exhibited in Figure~\ref{fig:test_1_p1c_u} where a continuous, smoother, projection of the actual discontinuous approximated solution had been plotted. As expected, blow-up is prevented for whatever choice of $c>0$ (even for small values like $c=10^{-3}$) if $\gamma$ satisfies the bound \ref{condgamma}, i.e., $1.5<\gamma\le 2$ (see Figure~\ref{fig:test_1_max-u_a}). However, in case that $\gamma$ does not comply with \ref{condgamma}, we may require a big enough value of $c$ to prevent the gathering phenomena. This value of $c$, apparently breaking down the tendency toward chemotactic collapse, increases as long as $\gamma$ moves away from the value $\gamma^*=1.5$ for which the equality holds in \ref{condgamma} (see Figures~\ref{fig:test_1_max-u_b} and \ref{fig:test_1_max-u_c}).

\subsubsection{Numerical simulations for the fully nonlinear attraction-repulsion model}
\begin{figure}
    \centering
    \begin{subfigure}[b]{0.3\textwidth}
    \centering
        \includegraphics[width=\textwidth]{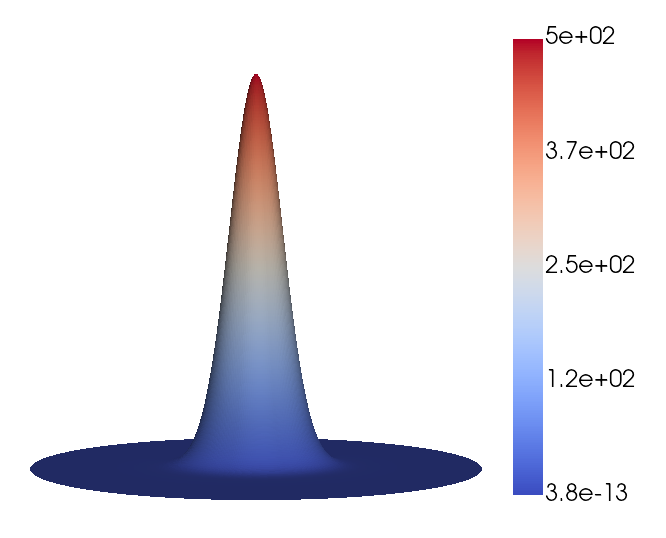}
        \caption{$t=0$.}
        \label{fig:test_2_p1c_u_a}
    \end{subfigure}
    \hfill
    \begin{subfigure}[b]{0.3\textwidth}
    \centering
        \includegraphics[width=\textwidth]{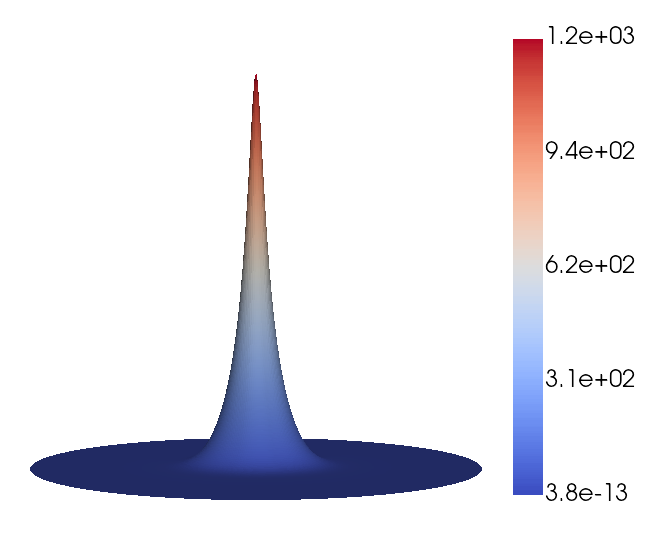}
        \caption{$t=5\cdot 10^{-5}$.}
        \label{fig:test_2_p1c_u_b}
    \end{subfigure}
    \hfill
    \begin{subfigure}[b]{0.3\textwidth}
    \centering
        \includegraphics[width=\textwidth]{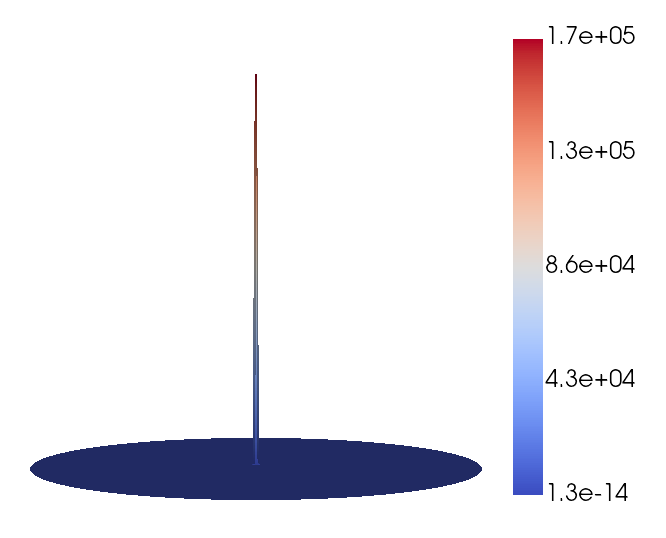}
        \caption{$t=3\cdot 10^{-4}$.}
        \label{fig:test_2_p1c_u_c}
    \end{subfigure}      
    \caption{Approximation of the solution $u$ at different time steps ($m_1=1$, $c=0$, $\alpha=1.5$).}
    \label{fig:test_2_p1c_u}
\end{figure}

Motivated by \cite{ColFraVig-StudApplMat}, let us analyze model \eqref{localproblem} under the assumption that 
\begin{equation}\label{cond}
m_2=m_3, m_1\in\R, \alpha>\beta \quad \text{and} \quad 
    \begin{cases}
m_2+\alpha > \max\bigl\{m_1+\frac{2}{n}k, k\bigr\} & \text{ if } m_1 \geq 0,\\
m_2+\alpha > \max\bigl\{\frac{2}{n}k, k\bigr\} & \text{ if } m_1 < 0.
\end{cases}
\end{equation}
In particular, we focus on the 2D case ($n=2$) in the domain $\Omega=\{(x,y)\colon x^2+y^2<1\}$. We set $\alpha=1.5$, $\xi=\chi=1$, and we define the following initial condition
$$
u_0(x,y)\coloneqq500e^{-35(x^2+y^2)},
$$
plotted in Figure~\ref{fig:test_2_p1c_u_a}. Also, we take a spatial and a temporal partition of sizes $h\approx 1.37\cdot 10^{-2}$ and $\Delta t=10^{-6}$.

\begin{figure}
    \centering
    \includegraphics[scale=0.5]{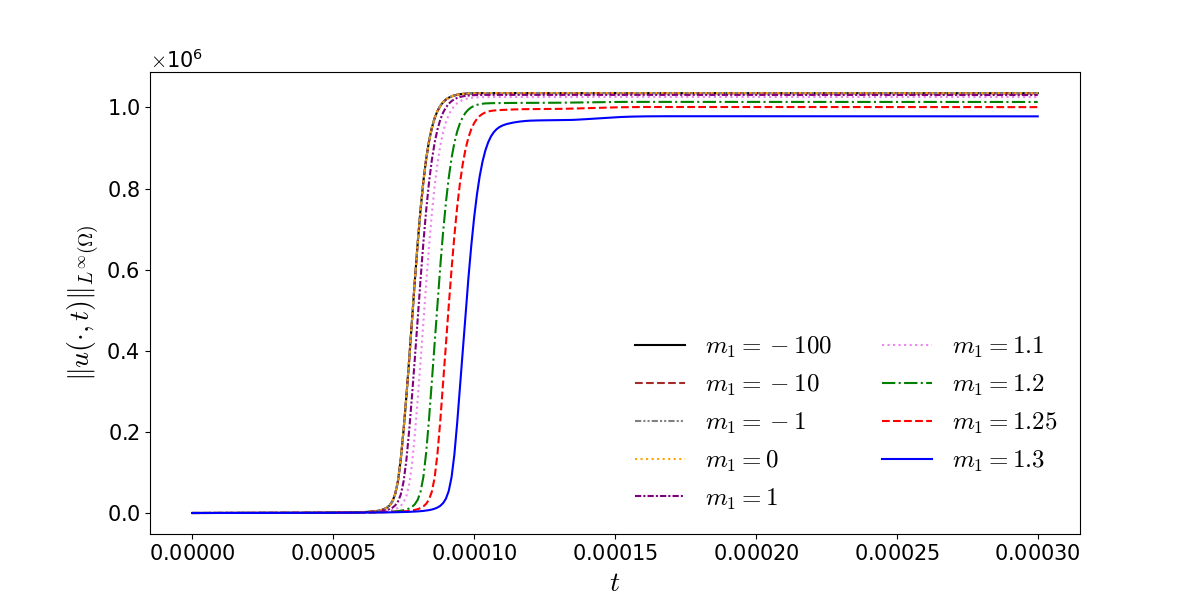}
    \caption{Maximum of the approximation of $u$ for different values of $m_1$ ($c=0$, $\alpha=1.5$).}
    \label{fig:test_1_max-u_m1}
\end{figure}

\begin{figure}
    \centering
    \includegraphics[scale=0.5]{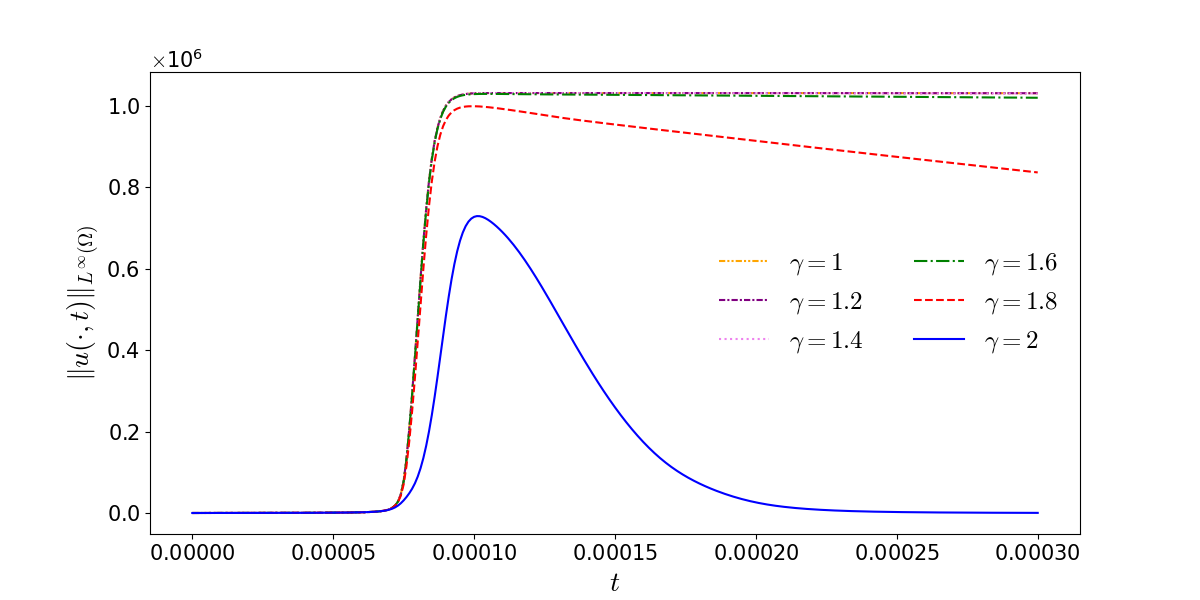}
    \caption{Maximum of the approximation of $u$ for different values of $\gamma$ ($m_1=1$, $c=10^{-3}$, $\alpha=1.5$).}
    \label{fig:test_1_max-u_gamma}
\end{figure}

First, in Figure~\ref{fig:test_2_p1c_u} we plot the evolution of $u$ over time without the damping gradient term, i.e., $c=0$. We observe that due to the choice of the parameters and the initial condition, it seems to occur a blow-up phenomenon (notable in Figures  \ref{fig:test_2_p1c_u_b} and \ref{fig:test_2_p1c_u_c}), in accordance with the results in \cite[Theorem 3]{ColFraVig-StudApplMat} under restriction \eqref{cond}.
In fact, if we use a nonlinear diffusion term moving the value of $m_1$ we can observe in Figure~\ref{fig:test_1_max-u_m1} that we still obtain a chemotactic collapse but with times of explosion. Consistently to the real phenomenon,  one may notice that the blow-up time increases (accordingly to Criterion \ref{CriterionBU}) with $m_1$, the diffusion coefficient working against coalescence effects, specially for large values.
Moreover, as in the previous example, if we introduce the damping gradient nonlinearities using a strictly positive small value for $c$, such as $c=10^{-3}$, the blow-up phenomenon is avoided if we stick to the bounds on $\gamma$ given in \ref{condgamma}, i.e., $1.67\approx5/3<\gamma\le2$. However, this singularity seems to remain for a small enough value of $c$ if $\gamma$ does not satisfy \ref{condgamma}; see Figure~\ref{fig:test_1_max-u_gamma}.

\subsubsection*{Acknowledgements}
The research of TL is supported by NNSF of P. R. China (Grant No. 61503171), CPSF (Grant No. 2015M582091), and NSF of Shandong Province (Grant No. ZR2016JL021).
The author DAS has been supported by UCA FPU contract UCA/REC14VPCT/2020 funded by Universidad de C\'adiz and by a Graduate Scholarship funded by the University of Tennessee at Chattanooga. DAS also acknowledges the travel grants to visit the Università di Cagliari funded by the Universidad de Cádiz and the Erasmus KA+131 program. The authors AC and GV are members of the Gruppo Nazionale per l’Analisi Matematica, la Probabilità e le loro Applicazioni (GNAMPA) of the Istituto Nazionale di Alta Matematica (INdAM) and  are partly supported by GNAMPA-INdAM Project \textit{Equazioni differenziali alle derivate parziali nella modellizzazione di fenomeni reali} (CUP--E53C22001930001) and by \textit{Analysis of PDEs in connection with real phenomena} (2021, Grant Number: F73C22001130007), funded by \href{https://www.fondazionedisardegna.it/}{Fondazione di Sardegna}. GV is also supported by MIUR (Italian Ministry of Education, University and Research) Prin 2022 \textit{Nonlinear differential problems with applications to real phenomena} (Grant Number: 2022ZXZTN2). AC is also supported by GNAMPA-INdAM Project \textit{Problemi non lineari di tipo stazionario ed evolutivo} (CUP--E53C23001670001).


\end{document}